\documentclass{amsart}
\usepackage{amssymb,amsmath,euscript}
\usepackage{hyperref}
 \theoremstyle{definition}
 
 \theoremstyle{remark}

\def\SL{\mathrm{SL}}
\def\GL{\mathrm{GL}}
\def\U{\mathrm{U}}
\def\OSp{\mathrm{OSp}}
\def\const{\mathrm{const}}

\def\diag{\mathrm{diag}}
\def\ov{\overline}

\def\cF{\mathcal F}

\DeclareMathOperator{\tr}{tr}
\DeclareMathOperator{\str}{str}

\newcommand{\A}{{\mathcal A}}

\newcommand{\C}{\ensuremath{\mathbb{C}}}
\newcommand{\R}[1]{{\mathbb R}^{#1}}

\newcommand{\RR}{\mathbb R}
\newcommand{\CC}{\mathbb C}
\newcommand{\Z}{{\mathbb Z_{2}}}
\newcommand{\ZZ}{{\mathbb Z}}
\newcommand{\fun}{C^{\infty}}
\DeclareMathOperator{\sdet}{sdet}

\DeclareMathOperator{\Ber}{Ber}
\DeclareMathOperator{\Diff}{Diff}

\newcommand{\der}[2]{{\frac{\partial {#1}}{\partial {#2}}}}

\newcommand{\dder}[3]{{\frac{\partial^2 {#1}}{\partial {#2}\partial {#3}}}}

\newcommand{\x}{{\xi}}
\newcommand{\h}{{\eta}}

\newcommand{\jtt}{{\tilde \jmath}}
\newcommand{\itt}{{\tilde \imath}}
\newcommand{\at}{{\tilde a}}

\begin{document}

\title[Felix Alexandrovich Berezin]
 {Felix Alexandrovich Berezin  and his work}

\author{Alexander Karabegov}
\address{%
Department of Mathematics\\
Abilene Christian University\\
Box 28012\\
Abilene, Texas 79699-8012\\
USA}
\email{alexander.karabegov@acu.edu}
%\thanks{17 (30) January 2012}

\author{Yuri Neretin}
\address{Institute for Theoretical and Experimental Physics\\
Bolshaya Cheremushkinskaya, 25 \\
117218 Moscow \\
Russia}
\email{neretin@mccme.ru}

\author{Theodore Voronov}
\address{School of Mathematics\\
University of Manchester\\
Oxford Road\\
Manchester, M13 9PL\\
United Kingdom}
\email{theodore.voronov@manchester.ac.uk}

\subjclass{Primary 01A70; Secondary 58A50, 58Z99, 22E66, 81S99}

\keywords{Representation theory, Laplace operators, second
  quantization, deformation quantization, symbols, supermanifolds, Lie
  supergroups}

\date{16 (29) February 2012}
\dedicatory{To the memory of F.~A.~Berezin (1931\,--\,1980)}
%%% ----------------------------------------------------------------------

\begin{abstract}
  This is a survey of Berezin's work focused on three topics:
  representation theory, general concept of quantization, and
  supermathematics.
\end{abstract}

%%% ----------------------------------------------------------------------
\maketitle
%%% ----------------------------------------------------------------------
%\tableofcontents
\section{Preface} This text has resulted from our participation in the
XXXth Workshop on Geometric Methods in Physics held in
Bia{\l}owie\.{z}a in summer 2011. Part of this conference was a
special Berezin Memorial Session: Representations, Quantization and
Supergeometry. F.~A.~Berezin, who died untimely in 1980 in a water
accident during a trip to Kolyma, would have been eighty in 2011.

This is an attempt to give a survey of Berezin's remarkable work and
its influence for today. Obviously, we could not cover
everything. This survey concentrates on three topics: representation
theory, quantization and supermathematics. Outside of its scope
remained, in particular, some physical works in which Berezin was
applying his approach to second quantization and his theory of
quantization. Also, we did not consider two important but somewhat
stand-alone topics of the latest period of Berezin's work devoted to
an interpretation of equations such as KdV from the viewpoint of
infinite-dimensional groups~\cite{berezin:kdvfa, berezin:kdvcmp}
(joint with A.~M.~Perelomov) and a method of computing characteristic
classes~\cite{berezin:charcl} (joint with V.~S.~Retakh).

For a sketch of Berezin's life and personality, we refer to a
brilliant text by R.~A.~Minlos~\cite{minlos}.

Sections~\ref{sec.laplace} and \ref{sec.secquant} below were written
by Yu. A. Neretin. Section~\ref{sec.quantiz} was written by
A.~V.~Karabegov. Section~\ref{sec.super} was written by
Th.~Th.~Voronov, who also proposed the general plan of the paper and
made the final editing.

\section{Laplace operators on semisimple Lie
  groups} \label{sec.laplace} The main scientific activity of
F.~A.~Berezin was related with mathematical physics, quantization,
infinite-dimensional analysis and infinite-dimensional groups, and
supermathematics.  But in 1950s he started in classical representation
theory (which at that time was new and not yet classical).

\subsection{Berezin's Ph.D. thesis: characters of complex semisimple
  Lie groups and classification of irreducible representations.}
Our first topic%
\footnote{This was not the first work of Berezin.  The paper \cite{BG}
  of Berezin and I.~M.~Gelfand (1956) on convolution hypergroups was
  one of the first attacks on the Horn problem; in particular they
  showed a link between eigenvalue inequalities and tensor products of
  irreducible representations of semisimple groups, see \cite{Klya},
  \cite{Ful}.} is the cluster of papers 1956-57:
announcements~\cite{L1}, \cite{L2}, \cite{bggn}, \cite{L3}, the main
text~\cite{Ber-main}, and an addition in~\cite{Ber-letter}. This work
has a substantial overlap with Harish-Chandra's papers of the same
years, see \cite{Harish2}.  F.~A.~Berezin in 1956 claimed that he
classified all irreducible representations of complex semisimple Lie
groups in Banach spaces.  We shall say a few words about this result
and the approach, which is interesting no less than the
classification.

The technology for construction representations of semisimple groups
(parabolic induction and principal series) was proposed by
I.~M.~Gelfand and M.~A.~Naimark in book~\cite{GN}. On the other hand,
Harish-Chandra~\cite{Harish1} in 1953 proved the `subquotient
theorem': each irreducible representation is a subquotient of a
representation of the principal (generally, non-unitary) series.

Consider a complex semisimple (or reductive) Lie group $G$, its
maximal compact subgroup $K$ and the symmetric space $G/K$. For
instance, consider $G=\GL(n,\C)$; then $K=\U(n)$ and $G/K$ is the
space of positive definite matrices of order $n$.  A {\it Laplace
  operator} is a $G$-invariant partial differential operator on $G/K$.
Let us restrict a Laplace operator to the space of $K$-invariant
functions (for instance, in the example above it is the space of
functions depending on eigenvalues of matrices). The {\it radial part
  of Laplace operator} is such a restriction.

Berezin described explicitly the radial parts of the Laplace operators
on $G/K$.  He showed that in appropriate coordinates%
\footnote{We also allow change $f(t)\mapsto \alpha (t)f(t)$.}  $t_1$,
\dots, $t_n$ on $K\!\setminus G/K$ each radial part has the form
 \begin{equation}
 p\Bigl(\frac{\partial}{\partial t_1}, \dots \frac{\partial}{\partial t_n}\Bigr)
 ,
 \label{eq:p}
 \end{equation}
 where $p$ is a symmetric (with respect to the Weyl group) polynomial.

 The first application was a proof of the formula for spherical
 functions on complex semisimple Lie groups from Gelfand and Naimark's
 book~\cite{GN}. One of the possible definitions of {\it spherical
   functions}: they are $K$-invariant functions on $G/K$ that are
 joint eigenfunctions for the Laplace operators. I.~M.~Gelfand and
 M.~A.~Naimark proved that for $G=\GL(n,\C)$ such functions can be
 written in the terms of the eigenvalues $e^{t_k}$ as
\begin{equation}
\Phi_\lambda(t)=\const(\lambda)\cdot
\frac{\det_{k,m} \{e^{\lambda_k t_m}\}}{\det_{k,m} \{e^{k t_m}\}}
\label{eq:spheric}
\end{equation}
as in the Weyl character formula for finite-dimensional
representations of $\GL(n,\C)$, but the exponents $\lambda_j$ are
complex. They wrote the same formula for other complex classical
groups, but it seems that their published calculation%
\footnote{It is very interesting, an integration in the Jacobi
  elliptic coordinates.}  can be applied only for $\GL(n,\C)$.
Berezin reduced the problem to a search of common eigenvalues of
operators (\ref{eq:p}) and solved it%
\footnote{The function $\alpha$ from a previous footnote %footnote on the previous page
  is the denominator of (\ref{eq:spheric}).}.

Next, consider Laplace operators on a complex semisimple Lie group
$G$, i.e.  differential operators invariant with respect to left and
right translations on $G$.  We can consider $G$ as a symmetric space,
it acts on itself by left and right translations, $g\mapsto h_1^{-1}g
h_2$, the stabilizer of the point $1\in G$ is the diagonal
$\diag(G)\subset G\times G$, i.e., we get the homogeneous space
$G\times G/ \diag(G)$. Note also that $G\times G/\diag(G)$ is the
complexification of the space $G/K$.  We again can consider the radial
parts of Laplace operators as the restrictions of Laplace operators to
the space of functions depending on eigenvalues $\lambda_j$.  Since
now eigenvalues are complex, the formula transforms to
 \begin{equation}
   p\left(\frac{\partial}{\partial t_1}, \dots \frac{\partial}{\partial t_n};
     \frac{\partial}{\partial \ov t_1}, \dots \frac{\partial}{\partial \ov t_n}\right)
   ,
 \label{eq:pp}
 \end{equation}
 where $p$ is separately symmetric with respect to holomorphic and
 anti-holo\-morphic partial derivatives\,\footnote{The eigenfunctions
   of (\ref{eq:pp}) are exponential and we have to symmetrize them
   because we need symmetric solutions.}.

 Recall that for infinite-dimensional representations $\rho$ the usual
 definition of the character $\chi(g)=\tr \rho(g)$ makes no sense,
 because an invertible operator has no trace. However, for {\it
   irreducible} representations of semisimple Lie groups and smooth
 functions $f$ with compact supports the operators $\rho(f)=\int
 f(g)\rho(g)$ are of trace class.  Therefore $f\mapsto \tr f(g)$ is a
 distribution on the group in the sense of L.~Schwartz.  This is the
 definition of the {\it character} of an irreducible representation.

 A character is invariant with respect to the conjugations $g\mapsto
 hgh^{-1}$. Also, it is easy to show that a character is an
 eigenfunction of all Laplace operators.  The radial parts of Laplace
 operators were evaluated, so we can look for characters as joint
 eigenfunctions of operators (\ref{eq:pp}). Algebraically the problem
 is similar to calculation of spherical functions and final formulas
 are also similar (but there are various additional analytic
 difficulties).

 For a generic eigenvalue, a symmetric solution is unique. It has the
 form
$$
\sum_{\sigma\in S_n} (-1)^\sigma e^{\sum_k(p_k t_{\sigma(k)}+ q_k\ov t_{\sigma(k)})}\,,
$$
for $G=\GL(n,\C)$, here $S_n$ is the symmetric group. This is the
character of a representation of the principal series.  For
`degenerate' cases there are finite subspaces of solutions.  Berezin
showed that all characters are linear combinations of the characters
of representations of principal series.  In the introduction
to~\cite{Ber-main}, he announced without proof a classification of all
irreducible representations.  The restriction of a representation of
the principal series to $K$ contains a unique subrepresentation with
the minimal possible highest weight%
\footnote{In 1966 D.~P.~Zhelobenko and M.~A.~Naimark~\cite{ZhN}
  announced the classification theorem in a stronger form. Later
  (1967--1973) D.~P.~Zhelobenko published a series of papers on
  complex semisimple Lie groups, e.g.~\cite{Zhe}, where he, in
  particular, presented a proof of this theorem (with a contribution
  of M.~Duflo).}.  We must choose a unique subquotient containing this
representation of $K$.

A formal proof of the classification of representations was not presented in~\cite{Ber-main},
but the theorem about characters and the classification theorem  are equivalent%
\footnote{It is not difficult to show that the distinct  subquotients have different
characters. The transition matrix between the characters of the principal series and
the characters of irreducible representations
is triangular with units on the diagonal.}.

\smallskip
{\small
Paper~\cite{Ber-main} was written in an enthusiastic style  and was not always careful. J.~M.~G.~Fell,  Harish-Chandra, A.~A.~Kirillov, and G.~M.~Mackey
formulated two critical arguments; Berezin responded in a separate paper~\cite{Ber-letter}.

Firstly, the original Berezin work contains a non-obvious and unproved lemma (on the correspondence between solutions of the systems of PDE in distributions on the group and the system of PDE in radial  coordinates). A proof  was a subject of the additional paper~\cite{Ber-letter}.

Secondly,  Berezin actually worked with irreducible representations whose $K$-spectra have finite multiplicities (i.e, the  irreducible Harish-Chandra modules). He formulated the final result as the 'classification of all {\it irreducible} representations in Banach spaces' and at this point   he
claimed  that    the equivalence of the two concepts had been proved by Harish-Chandra. But this  is not correct%
 \footnote{These two properties are not equivalent, see  Soergel's
 counterexample~\cite{Soe}.}.  He had to formulate the statement as the ``classification of all {\it completely
 irreducible}%
 \footnote{There are many versions of irreducibility for infinite-dimensional non-unitary
 representations. A representation is {\it completely irreducible} if
 the image of the group algebra is weakly dense in the algebra of all operators.}
  representations in Banach spaces", with the necessary implication  proved
 by R.~Godement~\cite{God} in 1952.

}

\smallskip
Recall that the stronger version of classification theorem was proved by Zhe\-lobenko near
1970. For real semisimple groups, the classification was announced by R.~Langlands in 1973 and
proofs were published by A.~Borel and N.~Wallach in 1980.

\subsection{Radial parts of Laplace operators.}

Spherical functions, the spherical transform, and the radial parts of Laplace operators appeared in representation
    theory in the 1950s.
Later  they became important in integrable systems.
On the other hand,
    they gave a new start
     for the theory of multivariable special functions (I.~G.~Macdonald, H.~Heckman, E.~Opdam,
     T.~Koornwinder, I.~Cherednik, and others.).

Consider a real semisimple Lie group $G$, its maximal compact subgroup $K$
  and the Riemannian symmetric space $G/K$. If the group $G$ is complex,
   then the spherical functions are
   elementary functions, as we have seen above.

But  for the simplest of the real groups, $G=\SL(2,\RR)$, the spherical functions are the Legendre functions.
In this case, the radial part of the Laplace operator is a hypergeometric differential operator (with some special values of the parameters).  General spherical functions are higher analogs of the Gauss hypergeometric functions.  Respectively, the radial parts of the Laplace operators are higher analogs of hypergeometric operators  (see expressions in~\cite{Sek} and~\cite{Heck}, Chapter 1).  The first  attack in this direction was made by   F.~A.~Berezin and F.~I.~Karpelevich~\cite{BK} in 1958.

Berezin and Karpelevich found a semi-elementary case, the pseudounitary group $G=\U(p,q)$. In this case the
radial parts of Laplace operators  are also symmetric expressions of the form
$$
r\bigl(L(x_1), \dots, L(x_p)\bigr)\,,
$$
but $L(x)$   is now a second order (hypergeometric) differential
operator,
$$
D:=x(x+1)\frac {d^2}{dx^2}+\bigl[(q-p+1)+(q-p)x\bigr] \frac {d}{dx}
+\frac 14(q-p+1)^2\,.
$$
They also evaluated the spherical functions on $\U(p,q)$ as eigenfunctions
 of the radial Laplace operators. In appropriate
coordinates the functions have the form
$$
\Phi_s(x)=\const\cdot\frac
{\det\limits_{k,j}
\bigl\{
 {}_2F_1\Bigl[\begin{array}{c}
\tfrac12(q-p+1)+is_j,\tfrac12(q-p+1)-is_j\\ q-p+1
\end{array};-x_k\Bigr]\bigr\}}
{\prod_{1\le k < l \le p}(s_k^2-s_l^2)
 \prod_{1\le k < l \le p}(x_k-x_l)}\,.
$$
Here $_2F_1[\dots]$ is the Gauss hypergeometric function, $x_1,\ldots, x_p$ are coordinates on the Cartan subgroup
of $\U(p,q)$, and $s_1, \ldots, s_p$ are parameters of spherical functions.

This paper was accepted by \emph{Doklady} in June 1957. Near that time Berezin's scientific interests
had changed and he left the classical representation theory%
\footnote{In 1976 paper~\cite{Ber-superupq} and the   five ITEP preprints of 1977 included in the
English version of \cite{Ber-super-book},   Berezin returned to the study of Laplace operators and considered the radial parts of Laplace operators for
Lie supergroups (see subsection~\ref{subsec.berezinlatest}). They are usual (non-super) partial differential operators. This topic is not well understood
up to now; A.~N.~Sergeev  and A.~P.~Veselov produced from this standpoint new operators of Calogero--Moser
type whose eigenfunctions are super-Jack functions (which also  are functions of even variables),
 see~\cite{Ser}. On an analog of the group case, see~\cite{Huck}.}$^,$%
 \footnote{In 1970s, Berezin made a work on the harmonic analysis
 in Hilbert spaces of holomorphic functions~\cite{Ber-hol1}, \cite{Ber2}, %\cite{Ber-hol2},
 %\cite{Ber-hol3};
 \cite{Ber10}
  for a discussion of this work and its continuations,
 see~\cite{UU}, \cite{Ner-Ber}, and \cite{Ner-gauss},  Chapter 7.}.

(The next step was done by M.~A.~Olshanetsky and A.~M.~Perelomov \cite{PO} in 1976;  see also \cite{PO2}. They wrote the radial part of the second order Laplace operator. Quite soon J.~Sekiguchi~\cite{Sek}
 obtained a  general formula for the groups $\GL$.)

\section{Method of second quantization}\label{sec.secquant}

Our next topic is the famous book ``The method of  second quantization"
\cite{Ber-second} (and  the announcements~\cite{Ber-second1, BMV, Ber-second2, Ber-second3}).
A more detailed   discussion of  the intellectual history of this work and its
influence  is in~\cite{Ner}.

\subsection{Prehistory.} It is known that at the end of 1950s Berezin started to
learn physics
and to participate in theoretical physics seminars in Moscow. He had to decide between numerous possible ways in this new world and his  choice was the problem about the automorphisms
of the canonical commutation and anticommutation relations formulated
in the book `Mathematical aspects of the quantum theory of fields' by
 K.~O.~Friedrichs~\cite{Fri} of 1953.

Let $P_1$, \dots, $P_n$, $Q_1$,\dots, $Q_n$ be self-adjoint operators in a Hilbert space
 satisfying the conditions
\begin{equation}
[P_k,P_l]=[Q_k,Q_l]=0,\qquad [P_k,Q_l]=i\delta_{k,l}
\label{eq:ccr}
\end{equation}
and without a common invariant subspace.
Such conditions are called the {\it canonical commutation relations}, abbreviation CCR.
According to the Stone--von Neumann theorem, such a system of operators is unique up to a unitary equivalence
(for a precise forms of the theorem, see, e.g.,~\cite{Ber-some-remarks}).
In fact, our Hilbert space can be identified with $L^2(\RR^n)$ and the operators with $x_k$,
 $i\frac \partial{\partial x_k}$, respectively.
 Now let
$g=\begin{pmatrix} \alpha&\beta\\ \gamma&\delta\end{pmatrix}$ be a symplectic $2n\times 2n$
matrix. Evidently, the operators
\begin{equation}
P'_k=\sum_l \alpha_{kl} P_l+ \sum_l \beta_{kl} Q_{kl},\qquad
Q'_k=\sum_l \gamma_{kl} P_l+ \sum_l \delta_{kl} Q_{kl}
\label{eq:transformation}
\end{equation}
satisfy the same relations (\ref{eq:ccr}). Therefore  there is a unitary operator $U=U(g)$ such that\,\footnote{Now the mapping $g\mapsto U(g)$ is called the \emph{Weil representation}, see A.~Weil's paper~\cite{Wei}, 1964.
The term is common and convenient, but historically
 it was a construction due to K.~O.~Friedrichs and I.~Segal.}
\begin{equation}
P'_k=U(g)P_kU(g)^{-1},\qquad Q'_k=U(g) Q_k U(g)^{-1}\,.
\label{eq:weil}
\end{equation}
By a version of the Schur Lemma, this operator is unique up to a scalar factor. It is easy to
see that
$$
U(g)U(h)=\lambda(g,h) U(g,h)\,,
$$
where $\lambda(\cdot,\cdot)$ is a complex scalar. Apparently, Friedrichs decided that there was nothing
to discuss here and asked what would happen if the number $n$ of the operators were $\infty$. He showed
that there are many nonequivalent representations of CCR besides the well-known Fock representation. Next, Friedrichs asked, for which symplectic matrices the system of operators $P'_k$, $Q_k'$ are equivalent
 to $P_k$, $Q_k$. He formulated a correct conjecture and tried to find explicit formulas for $U(g)$.

%%%%%%%%%%%%%%%%%%%%%%%%%%%%%%%%%%%%%%%%%%%%%%%%%%%%%%%%%%%%%%%%%%%%%%%%%%%%%%%%%%%%%%%%%%%%%%

\subsection{Operators and divergences.} Consider the usual Fourier transform $\cF$ in $L^2(\RR)$,
$\widehat f(\xi)=\int e^{ix\xi} f(x)\,dx$.
Its definition is not completely straightforward, since the integral can be divergent,
and some regularization dance  is necessary. If we want to find
$\cF^2$, we must calculate the kernel
$$
K(x,y)=
\int e^{iy\xi} e^{ix \xi} \,d\xi
$$
Since we know the answer, we can believe that it is obvious. In any case, the integral diverges...

These difficulties  are usual for the work with  integral operators in $L^2(\RR^n)$.
Field theory requires  functions of infinite number of variables and passing to
the limit   $n\to\infty$ only multiplies the problems.
Berezin noticed that in the space $F_n$ of entire functions on $\C^n$
with the inner product
$$
\langle f,g\rangle=\frac1{(2\pi)^n}\int f(z)\ov{g(z)} e^{-|z|^2}\,d\Re(z)\,d\Im(z)
$$
we can realize
our operators as
$$
P_k=\frac 1{\sqrt 2}\Bigl(z_k+\frac \partial{\partial z_k}\Bigr),
\qquad Q_k= \frac 1{\sqrt 2 i}\Bigl(z_k-\frac \partial{\partial z_k}\Bigr)\,.
$$
Therefore this space can be identified with $L^2(\RR^n)$.
Berezin observed that in the space $F_n$ any bounded  operator is an integral operator of the form
$$
Af(z)=\int_{\C^n} K(z,\ov u) f(u) e^{-|u|^2} du\,d\ov u
$$
and the integral is convergent. Also the kernel of a product of integral operators
is defined by  a convergent integral. Next, Berezin showed that this `holomorphic model'  perfectly survives
as $n\to\infty$ (only the case $n=\infty$ is discussed in book~\cite{Ber-second}, Berezin uses the term
`generating functional' for the function assigned to an operator). In particular, we can work with bounded operators without any divergences.
%Next, Berezin propose a parallel (equivalent) language of Wick symbols. Recall that expressions
%f operators in the form $A=p(x,\frac{\partial}{\partial x}$ were widely discussed in 1930s.%
%Expressions of the form $A=p(z,\frac{\partial}{\partial z}$ seems as same from algebraic
%point of view but they are essentially better from the

Certainly, we need also unbounded operators, where divergent expressions have to appear. But,
again, the `level of divergences' is minimal.

In parallel, Berezin proposed an almost equivalent formalism of Wick symbols.
Algebraically, they looks similar to the well-known since 1930s expressions of operators
as $A=p(x,\frac{\partial}{\partial x})$, where all $x$'a are at the left and all $\der{}{x}$'s are at the right. But only few operators can be written in this form
if we understand `functions' literally. In contrast, we can express an operator as
$A=p(z,\frac{\partial}{\partial z})$ more or less always.

\subsection{Weil representation.} Using this operator formalism, Berezin wrote explicit formulas
for the  operators $U(g)$. He interpreted the conditions (\ref{eq:weil}) as a first order
system of PDEs for the
kernels $K$ of $U(\cdot)$, solved the equations and got the expressions of the form
\begin{equation}
K(z,\ov u)=\exp\bigl\{ S(z, \ov u)\bigr\}\,,
\label{eq:Q1}
\end{equation}
where $S$ is an explicit quadratic form.
Thus we obtain a projective representation of an infinite-dimensional symplectic group by integral operators acting in
the space of functions of infinite number of variables. We also can replace $\infty\mapsto n$
and obtain a  construction that was completely new in that time.
%which now is called Weil representation%
%\footnote{Historically, it is Friedrichs--Segal representation, see above and below.}

In particular, Berezin proved the Friedrichs conjecture about the domain of definition
of this representation.

%%%%%%%%%%%%%%%%%%%%%%%%%%%%%%%%%%%%%%%%%%%%%%%%%%%%%%%%%%%%%%%%%%%%%%%%%%%%%%%%%%%%%

\subsection{Fermionic Fock space.} For us a fermionic Fock space is a space of functions
of anticommuting variables. This  idea, now common,  originated from Berezin's book~\cite{Ber-second}. Berezin also
found that there is a natural integral over anticommuting
variables (\cite{Ber-second1}). We say more about that in section~\ref{sec.super}.
Berezin showed that an operator in the fermionic Fock space is determined by a function
(the `generating functional') depending on a double collection
of anticommuting variables and that  it is convenient to express operators in a fermionic Fock
space as  integral operators, with respect to that peculiar integral.

In~\cite{Fri}, Friedrichs also formulated a problem about the \emph{canonical anticommutation relations}
(abbreviation CAR)
\begin{equation}
\{P_k,P_l\}=\{Q_k,Q_l\}=0,\qquad \{P_k,Q_l\}=i\delta_{k,l}
\label{eq:acr}
\end{equation}
and their symmetries (\ref{eq:transformation}). Now the matrix
$g=\begin{pmatrix} \alpha&\beta\\ \gamma&\delta\end{pmatrix}$ is orthogonal.
Berezin solved this problem as well and wrote a formula
for the kernels of $U(\cdot)$,
\begin{equation}
K(\xi,\eta)=\exp\{S(\xi,\eta)\}\,,
\label{eq:Q2}
\end{equation}
where $S$ is an explicit quadratic expression.
Note that the formulas for $S$ in (\ref{eq:Q1}) and (\ref{eq:Q2})
are similar.

In fact, the both theorems  are results in the representation theory of infinite-dimensional Lie groups.
Berezin's book  can be regarded as a mathematization of   field theory.
However,  it  was also (Chapters 2 and 3) the first book on infinite-dimensional
 groups and the start of this theory.  For a more detailed discussion, see~\cite{Ner}.

\subsection{History and references.}\label{subsec.secquant.history}
Main Berezin's results with outlined proofs were announced in   \emph{Doklady} paper~\cite{Ber-second1},
of March 1961 (accepted in November 1960). The text was written in the telegraphic style usual for   \emph{Doklady} of that time: the  allowed four pages  were all used up to one line.  In September 1962, Berezin submitted a large paper
to \emph{Uspehi} (that is, \emph{Russian Mathematical Surveys}). The paper was rejected. In the following years,
Berezin published more short announcements:~\cite{Ber-second2}, \cite{Ber-second3} and \cite{BMV}. In 1965\,%
\footnote{The English version appeared in 1966.} the book ``The method of second quantization'' was published, addressed to physicists%
\footnote{In spite of its physical language, the book is a rigorous, maybe not detailed, mathematical text.}.

Friedrichs's questions also attracted Irving Segal, who had  worked in mathematical
field theory since the beginning of 1950s.  (In particular, Segal introduced a model of the Fock space
as $L^2$ on a Gaussian measure~\cite{Seg2}, 1956;  later J.~Feldman~\cite{Fel}, 1959, constructed the action of an infinite-dimensional $\GL$ on that space.) In 1959, Segal obtained explicit formulas
for the `Weil representation' for finite $n$ in the space $L^2$, \cite{Seg}. In 1961, he proposed a holomorphic model
for the boson Fock space (this was also done by V.~Bargmann~\cite{Bar} in the same year).
In 1962, D.~Shale~\cite{Sha} published the solution of the Friedrichs problem for CCR,
and in 1965, D.~Shale and R.~W.~Stinespring published their solution
for CAR%
\footnote{Note that the famous mathematicians D.~Shale, R.~W.~Stenespring and J.~Feldman were all I.~Segal's students.} \cite{ShSt}.

However these papers did not cover Berezin's results. His book and Berezin himself immediately became famous.

\subsection{Berezin's book in physics.}
Besides the formal results concerning CCR and CAR,
the interest of physicists to this text had two additional reasons.

First, the new operator  formalism  (both bosonic and fermionic) was very convenient.
It became easier to write formulas and to calculate.

The second reason was  the mysterious  parallelism between the bosonic and fermionic spaces
which was emphasized in the book. For Berezin himself this was the starting point  of his work
leading to the creation of supermathematics (see section~\ref{sec.super}).

\section{Berezin's general concept of quantization}\label{sec.quantiz}

One of the main directions of Berezin's research was mathematical formulation of the concept of quantization as a deformation of a classical mechanical system. In \cite{Ber5} Berezin interpreted the universal enveloping algebra of a Lie algebra as a quantization of the Poisson algebra of polynomial functions on the dual of the Lie algebra. In \cite{Ber1} and \cite{Ber3} Berezin introduced a general concept of quantization based upon algebras of operator symbols depending on a small parameter.

Quoting from \cite{berezin:obituary-umn}, ``according to the main idea of these works, quantization has the following precise mathematical meaning: the algebra of quantum observables is a deformation of the algebra of classical observables, so that the Planck constant plays the role of the deformation parameter and the direction of deformation (the first derivative in the parameter at zero) is given by the Poisson bracket''

In \cite{Ber2}, Berezin studied quantization of complex symmetric spaces.  The operator symbols used in quantization were introduced and studied in \cite{Ber7}, \cite{Ber4}, and \cite{Ber9}.  In \cite{Ber10}, Berezin obtained the spectral decomposition of the operator connecting covariant and contravariant symbols on classical complex symmetric spaces, now called the Berezin transform.
In \cite{Ber6} and \cite{Ber8}, Berezin constructed finite approximations of Feynman path integrals with the use of operator symbols. See also Berezin and M.~A.~Shubin~\cite{berezin:ishubin1} and their joint book ``The Schr\"{o}dinger Equation''~\cite{berezin:ishubin}, which Shubin prepared for publication after Berezin's death.

Let us consider these works in greater detail.

\subsection{Poisson bracket and quantization on the dual of a Lie algebra.}
In the fundamental paper \cite{Ber5} Berezin constructed an integral transform $\delta$ from generalized functions on a neighborhood of the identity of a Lie group $G$ to functions on the dual $\tilde {\mathfrak G}$ of its Lie algebra ${\mathfrak G}$ and expressed the symmetrization mapping $\Lambda$ from the symmetric algebra $S$ of ${\mathfrak G}$ to the universal enveloping algebra $\hat S$ of ${\mathfrak G}$ through the map\-ping~ $\delta$. Let $\{\hat x_p\}$ be a basis in ${\mathfrak G}$, $\{t_p\}$ the corresponding coordinates on ${\mathfrak G}$, $\{y_p\}$ the dual coordinates on $\tilde {\mathfrak G}$. The mapping $\delta$ is defined as follows:
\[
     \delta: s(g) \mapsto\, \stackrel{\delta}{s}(y) = \int e^{-ity} s(g(t))\rho(t) dt, \quad ty = \sum_p t_p y_p,
\]
where $t \mapsto g(t)$ is the exponential mapping and $\rho(t)$ is the density of the right-invariant measure on $G$ in the canonical coordinates $\{t_p\}$. The symmetric algebra $S$ of ${\mathfrak G}$ is identified with the space of polynomials on $\tilde {\mathfrak G}$ and $\Lambda$ maps $y_p$ to $\hat y_p = - i \hat x_p$. The generalized functions supported at the identity of $G$ form an algebra with respect to the convolution. This algebra is naturally identified with the universal enveloping algebra $\hat S$. Berezin proved that under this identification the inverse mapping $\Lambda^{-1}$ is given by the mapping $\delta$. The mapping $\delta$ allows to transfer the convolution of generalized functions supported on a small neighborhood $U$ of the identity of the group $G$ to an operation on functions on the dual $\tilde {\mathfrak G}$ of the Lie algebra ${\mathfrak G}$. Berezin gave an integral formula for this operation. Given generalized functions $s_1,\ s_2$ supported on $U$ and their convolution $s$, set $\sigma_1 = \stackrel{\delta}{s_1},\sigma_2 = \stackrel{\delta}{s_2}$ and $\sigma = \stackrel{\delta}{s}$. Then
\[
   \sigma(y)=\int K_U(y\vert y_1,y_2)\sigma_1(y_1)  \sigma_2(y_2) dy_1 dy_2,
\]
where
\[
          K_U(y\vert y_1,y_2) =  \frac{1}{(2\pi)^{2n}} \int_{g(t_1)\in U,\, g(t_2)\in U}
                  {e}^{-iy\log(g(t_1)g(t_2)) + iy_1 t_1 + iy_2t_2} dt_1 dt_2.
\]
Berezin noted that this integral formula can be extended to the space $S$ of polynomials on $\tilde {\mathfrak G}$ and the resulting algebra is isomorphic to the universal enveloping algebra $\hat S$ of ${\mathfrak G}$. Moreover, the leading term of the commutator of polynomials leads to a natural Poisson bracket on the dual of the Lie algebra ${\mathfrak G}$. For arbitrary smooth functions on $\tilde {\mathfrak G}$ it is possible to write
\begin{equation*}
    \{f_1,f_2\} = \sum  C_{ij}^k y_k \der{f_1}{y_i}\der{f_2}{y_j}\,,
\end{equation*}
where
 $C_{ij}^k$
 are the structure constants of the Lie algebra ${\mathfrak G}$. About the same time, this Poisson bracket on $\tilde {\mathfrak G}$ (in the form of a symplectic structure on the coadjoint orbits) was discovered in the orbit method. Therefore it became known as the \emph{Berezin--Kirillov } or \emph{Berezin--Kirillov--Kostant bracket} . (Later Alan Weinstein found out that the bracket had been known already to S.~Lie, so the name `Lie--Poisson bracket' became more standard.)

Thus the universal enveloping algebra of ${\mathfrak G}$ can be interpreted as a quantization of the corresponding Poisson algebra on~$\tilde {\mathfrak G}$ consisting of the polynomial functions endowed with  the Berezin--Kirillov--Kostant Poisson bracket.

\smallskip

{\small
In fact, by rescaling this operation by a formal parameter $\hbar$, one can obtain from   above Berezin's formula the following integral formula for what is now known as the `Baker--Campbell--Hausdorff star product' on $\tilde {\mathfrak G}$\,:
\begin{multline} \label{eq.bch}
(f_{1}*f_{2})(y) = \\
\frac{1}{(2\pi\hbar)^{n}}\iint dy_{1}\,dt_{1}\,dy_{2}\,dt_{2} \;
f_{1}(y_{1})\,f_{2}(y_{2})\;
e^{-\frac{i}{\hbar}\left(\langle t_{1},y_{1}\rangle +
\langle t_{2},y_{2}\rangle - \langle H(t_{1}, t_{2}), y\rangle\right)}\,.
\end{multline}
Here $H(t_1,t_2)$ is the formal BCH power series on ${\mathfrak G}$ and the integration extends over $\tilde {\mathfrak G}\times\mathfrak G\times\tilde {\mathfrak G}\times\mathfrak G$. The functions $f_1$ and $f_2$ can be arbitrary smooth functions on~$\tilde {\mathfrak G}$ due to the presence of a formal parameter $\hbar$\,\footnote{We obtained formula~\eqref{eq.bch} around 1998 (Th.V., unpublished) and then realized that it can be deduced from Berezin~\cite{Ber5}.}.

}

\subsection{General concept of quantization as deformation.}
In \cite{Ber1} and \cite{Ber3} Berezin gave a general definition of quantization of a Poisson manifold
$(M, \{\cdot,\cdot\})$ as an algebra $({\mathfrak A},\ast)$ of sections of a field of noncommutative algebras $(\A_h,\ast_h)$ parameterized by the elements $h$ of a set $E$ of positive numbers that has zero as an accumulation point. The Correspondence Principle for this quantization is expressed in terms of a homomorphism
\begin{equation}\label{E:homphi}
\varphi_0: {\mathfrak A} \to C^\infty(M)
\end{equation}
such that for $f,g \in {\mathfrak A}$,
\[
    \varphi_0\left(\frac{1}{h}(f\ast g - g \ast f)\right) = i\{\varphi_0(f),\varphi_0(g)\}.
\]
Then he considered a special case when $\A_h \subset C^\infty(M)$, the elements of ${\mathfrak A}$ are functions $f(h,x)$ on $E \times M$, and
\[
    \varphi_0(f) = \lim_{h\to 0} f(h,x).
\]

\subsection{Berezin's quantization using symbols.}
Berezin studied a number of examples of such special quantizations where $\A_h$ for a fixed $h$ is an algebra of symbols of operators in a Hilbert space. To this end, Berezin introduced \emph{covariant} and \emph{contravariant  symbols} related to an overcomplete family of vectors in a reproducing kernel space. Namely, consider a Hilbert space $H$ and a set $M$ with measure $d\alpha$ whose elements parameterize a system of vectors $\{e_\alpha\}$ in $H$. Let $P_\alpha$ be the orthogonal projection operator onto $e_\alpha$ and
\[
    d\mu(\alpha) = ||e_\alpha||^2 d\alpha
\]
be another measure on $M$. The vectors $\{e_\alpha\}$ form an overcomplete family in $H$ if
\[
     \int P_\alpha\, d\mu(\alpha) = E
\]
is the identity operator in $H$. Then $H$ is isometrically embedded into $L^2(M, d\alpha)$ by the mapping $H \ni f \mapsto \langle f, e_\alpha \rangle$.
The projectors $P_\alpha$ are used to define covariant and contravariant symbols of operators in $H$ as follows. The covariant symbol of an operator $\hat A$ is the function
\[
A(\alpha) = \tr \hat AP_\alpha
\]
on $M$. A function $\stackrel{\circ}{\text{A}}(\alpha)$ on $M$ is a contravariant symbol of $\hat A$ if
\[
         \hat A = \int P_\alpha\, \stackrel{\circ}{\text{A}}(\alpha) d\mu(\alpha).
\]
The measure $\mu$ defines a trace functional on appropriate classes of covariant and contravariants symbols that agrees with the operator trace (see \cite{Ber4}),
\[
    \tr \hat A = \int A\, d\mu = \int \stackrel{\circ}{\text{A}}\, d\mu.
\]
The covariant and contravariant symbols $A$ and $\stackrel{\circ}{\text{A}}$ of the same operator $\hat A$ are connected via the Berezin transform $I$,
\[
     A(\alpha) = (I\stackrel{\circ}{\text{A}})(\alpha) = \int \tr \left(P_\alpha P_\beta\right) \stackrel{\circ}{\text{A}}(\beta) d\mu(\beta).
\]

An overcomplete system of vectors $\{e_\alpha\}$ in $H$ may admit a symmetry group $G$ that acts upon $H$ by a unitary representation $g \mapsto U_g$ and upon $M$ by transformations preserving the equivalence class of the measure $d\alpha$ so that
\[
     U_g e_\alpha = s(\alpha, g) e_{g\alpha},
\]
where $s: M \times G \to \C$ is a measurable cocycle satisfying
\[
      \frac{dg\alpha}{d\alpha} = |s(\alpha,g)|^2.
\]
Then $U_g P_\alpha U_g^{-1} = P_{g\alpha}$, the measure $d\mu$ is $G$-invariant, the symbol mappings
\[
     \hat A \mapsto \tr \hat A P_\alpha \mbox{ and } \stackrel{\circ}{\text{A}}(\alpha) \mapsto  \int P_\alpha\, \stackrel{\circ}{\text{A}}(\alpha) d\mu(\alpha)
\]
are $G$-equivariant and the Berezin transform $I$ is $G$-invariant.

Berezin studied spectral properties of covariant and contravariant symbols in \cite{Ber4} and then used algebras of covariant symbols to define a quantization of a special class of K\"ahler manifolds in \cite{Ber1} using the saddle-point method. He started with a K\"ahler manifold $M$ of complex dimension $m$ with a K\"ahler form $\omega$ and the Liouville measure $\omega^m$. He assumed that there exists a global K\"ahler potential $\Phi$ of the form $\omega$ and introduced an $h$-parameterized family of measures
\[
     d\alpha_h = {\rm e}^{-\frac{1}{h}\Phi} \omega^m
\]
on $M$. Then he considered the Hilbert space $H_h$ of holomorphic functions on $M$ square integrable with respect to the measure $d\alpha_h$. The Bergman reproducing kernel of $H_h$ defines an overcomplete system of vectors $\{e^{(h)}_\alpha\}$  in $H_h$. In order to prove the Correspondence Principle, Berezin imposed a severe assumption on $M$ that
\[
      e^{(h)}_\alpha (z) = c(h) {\rm e}^{\frac{1}{h}\Phi(z, \bar w)}
\]
for $\alpha = (w,\bar w) \in M$ and some constant $c(h)$. This assumption is satisfied on K\"ahler manifolds with a transitive symmetry group which allowed Berezin to quantize complex symmetric spaces (see \cite{Ber2}).

\subsection{Influence of Berezin's work.}

In the following decades Berezin's work on quantization attracted a lot of attention. His results were expanded and generalized by many mathematicians and mathematical physicists in two major directions. First, Berezin's definition of quantization in the special case when $\A_h \subset C^\infty(M)$ was extended to incorporate deformation quantization of Flato et al. \cite{BFFLS} as a formal asymptotic expansion in $h$ of the product $\ast_h$ in the algebra $\A_h$. In the general case this can be achieved by extending the homomorphism (\ref{E:homphi}) in the general definition of quantization to a homomorphism
\[
   \varphi = \varphi_0 + \nu \varphi_1 + \ldots: {\mathfrak A} \to C^\infty(M)[[\nu]]
\]
to the star-algebra of some formal deformation quantization on the Poisson manifold $(M, \{\cdot,\cdot\})$ such that $\varphi(hf) = \nu \varphi(f)$\footnote{Here $\nu$ is a formal parameter.}. Examples of such quantizations of K\"ahler manifolds were first given in \cite{Mor} and
 \cite{CGR1, CGR2}. The second direction was to remove the restrictions on the K\"ahler manifold in Berezin's quantization. Based on the microlocal technique developed by Boutet~de~Monvel and Guliiemin in~\cite{BMG}, it was shown in \cite{BMS} that Berezin--Toeplitz quantization\footnote{\emph{Berezin--Toeplitz quantization} is defined in terms of operators with given contravariant symbols. Such operators are generalizations of Toeplitz operators.}  on general compact K\"ahler manifolds satisfies an analog of the Correspondence Principle. Then in \cite{Sch} the existence of the corresponding Berezin--Toeplitz star product was established.  In \cite{CMP1} all star products ``with separation of variables" on an arbitrary K\"ahler manifold were classified and in \cite{KSch} the Berezin--Toeplitz star product was completely identified in terms of this classification. In \cite{E} M.~Engli\v{s} showed the existence of Berezin star-product on a quite general class of inhomogeneous complex domains. Berezin--Toeplitz quantization was recently studied by microlocal methods developed in \cite{MM} and \cite{LCh}.  Applications of Berezin--Toeplitz quantization in the topological quantum field theory were given in \cite{JEA}.

Much work has been done to generalize Berezin's quantization on K\"ahler manifolds to other spaces. Berezin's first doctoral student Vladimir Molchanov developed harmonic analysis and quantization on para-Hermitian symmetric spaces (see \cite{Mol1}, \cite{Mol2}). Berezin's quantization on quantum Cartan domains was considered in \cite{SSV}. Berezin's quantization was generalized to supermanifolds in \cite{BKLR}, \cite {GAN}. In the framework of this publication it is impossible to give a comprehensive survey of the growing body of papers building upon Berezin's work on quantization and many important papers are inevitably left out.

\section{Supermathematics}\label{sec.super}

\subsection{Introductory remarks.}
Without doubt, Berezin is \emph{the} creator of supermathematics, though it was not him who introduced the name. (More about the origin of the name in~\ref{subsec.supersym}.) In   hindsight, it is possible to trace the origins of what became supermathematics in various areas of pure mathematics and theoretical physics, but it only due to Berezin's vision and  his conscious effort  that these previously disjoint pieces became parts of a great unified picture together with a lot of new mathematics discovered by Berezin himself and by those who followed him. Speaking about Berezin's work in supermathematics, it is worth pointing out that it was interrupted  by his untimely death when  supermathematics was still in the early stages of its  development; therefore, the loss caused by Berezin's  sudden departure was   greater for supermathematics than for other areas of his work.

Berezin's publications related to supermathematics can be divided into two groups corresponding to the two periods: the gestation period (1961--1975) and the `super' period\,\footnote{The `super' period is marked by the emergence of the names such as supermanifold, superalgebra, etc. As the borderline we may take the discovery of supersymmetric physical models followed by the introduction of the notion of a supermanifold in mathematics. This division is partly conventional.} (1975--1980).

We can formulate the main idea of supermathematics as follows. The systematic consideration of $\Z$-graded objects such as Abelian groups, vector spaces, algebras and modules with the corresponding sign convention (``Koszul's sign rule''\,\footnote{From an abstract viewpoint, the sign rule used in supermathematics is a very special example of a ``commutativity constraint'' or ``braiding'' in  tensor categories.}) allows to construct a natural extension of the `usual' linear algebra including generalizations of commutative algebras and Lie algebras. This goes further to the extension of differential and integral calculus of many variables and, geometrically,  to the  extensions of the notions of differentiable manifold, Lie group, algebraic variety (or scheme) and   algebraic group.

Two things should be said.

Firstly, Berezin came to his program of supermathematics (without such a name, which appeared later) motivated by physics, more precisely, by his studies of the formalism of second quantization, which lies in the foundations of quantum field theory. The influence of physics was also decisive for the passage of supermathematics from its gestation stage to the  modern  stage.

Secondly, the `supermathematical' generalization of the usual notions is \emph{not arbitrary}, but indeed reflects the nature of things:  the `superanalogs' of various objects fit together in the same way as their prototypes do (but may also show non-trivial new phenomena). Moreover, this generalization is \emph{rigid and unique}. There are no known further generalizations based on other gradings or more complicated commutativity constraints. That is, there are quantum groups and quantum spaces; however, they  are isolated examples unified philosophically but not by a  general theory such as  a (non-existing) `quantum' or `braided' geometry, although these terms are sometimes applied.

\subsection{Analysis on a Grassmann algebra.}
As already said, Berezin's program of supermathematics has its roots in his book~\cite{Ber-second} and the related articles~\cite{Ber-second1, BMV, Ber-second2, Ber-second3}. (See the historical remarks  in~\ref{subsec.secquant.history}.) In order to construct a `calculus of functionals' for the Fermi fields that can be parallel to the functional calculus used  for describing the Bose fields, Berezin introduced differentiation and integration on a Grassmann algebra. He did it first for the Grassmann algebras with finite numbers of generators and then  extended the results to the infinite-dimensional `functional' case that he needed for his problem.
This calculus allowed Berezin to obtain a `functional realization' of  the fermionic Fock space similar to the realization of  the bosonic Fock space by holomorphic or antiholomorphic functions (or functionals) and to construct the spinor representation of the canonical transformations in the fermionic case (i.e., the spinor representations of certain infinite-dimensional versions of the orthogonal group). These representations were   discussed in Section~\ref{sec.secquant} and  we shall not repeat it here.

A striking feature of Berezin's calculus on a Grassmann algebra was, as he noted in~\cite{Ber-second},   that \emph{``differently from the usual rule for a change of variables, the independent variables and their differentials transform by reciprocal matrices''}. Here Berezin refers to his formulas
\begin{equation*}
    \int\!x\,dx=1\,, \quad \int\!dx=0
\end{equation*}
(where $x$ is a Grassmann generator), which imply $d(ax)=\frac{1}{a}\,dx$.  The `differential $dx$'    is   the quantity appearing under the integral sign, which we would now call the  Berezin volume element and denote by $Dx$ to distinguish it from the genuine differential of the variable $x$. We can see here the origins of  superdeterminant (now called Berezinian), which was discovered   some  years later\,\footnote{When  only odd variables are present, the Berezinian of their linear homogeneous transformation reduces simply to the inverse of the determinant.}.    Another remarkable fact noted and used by Berezin in~\cite{Ber-second} was the appearance of the (square root of the) determinant of the quadratic form for a `fermionic Gaussian integral' in the numerator, not in the denominator:
\begin{equation*}
    \int e^{\sum a_{ik}x_ix_k}\,dx_n\ldots dx_1= (\det \|2a_{ik}\|)^{1/2}\,,
\end{equation*}
where $a_{ik}=-a_{ki}$,   in a sharp contrast with the familiar (``bosonic'') case (equation~3.16 in~\cite{Ber-second}).

Integration on a Grassmann algebra introduced by Berezin was soon applied by Faddeev and Popov in their famous work on quantization of the Yang--Mills field: they expressed the Jacobian factor arising from the separation of the gauge degrees of freedom by a fermionic Gaussian integral over `ghosts' (the ``Faddeev--Popov ghosts'') and thus they were able to deduce  the Feynman rules including  ghosts  as following from a local Lagrangian field theory. By contrast, Bryce DeWitt, who obtained  close results   at the same time, did not know Berezin's integration and because of that failed to arrive
to such a natural formulation\,\footnote{DeWitt explicitly admits that in the preface to the Russian translation (1985) of his influential book `The dynamic theory of groups and fields'.}.

At this point it makes sense to discuss the question about Berezin's predecessors. It is sometimes claimed that the use of anticommuting variables for the classical description of fermions was familiar to quantum physicists since 1950s (and hence Berezin did not introduce anything particularly new). Typically Schwinger's name is mentioned in this regard. In reality, the ideas of Schwinger and his disciples such as DeWitt about anticommuting variables were quite vague and did not go any further than the introduction of `left' and `right' derivatives with respect to generators of a Grassmann algebra (see, e.g., \cite{schwinger:dynamprinc53, schwinger:grassmann62}). Partial differentiation with respect to exterior generators had been already known to \'{E}lie Cartan in connection with his method of \emph{rep\`{e}re mobile} (physicists were probably unaware of that). The novelty of Berezin's work in comparison with earlier and simultaneous works by physicists was in the mathematical clarity  and power  with which he developed  the analogy between usual functions and elements of a Grassmann algebra, but the main  new feature was   integration over Grassmann generators with its striking properties. There is a saying\,\footnote{I learned it from A.~A.~Kirillov.} that ``derivatives are  algebra; analysis begins with   integrals''. It was Berezin who made this decisive step. It took some time for this achievement to be absorbed by the physical community: the example of DeWitt
is a clear evidence.

\smallskip
{\small
One person who can be counted as a true predecessor of Berezin, is the  British physicist J.~L.~Martin. In two  papers~\cite{martin:classical, martin:feynman} of 1959, Martin introduced the notion of a general Hamiltonian system on a Poisson manifold (in the modern terminology) and suggested to extend it to more general algebras; in particular, he showed  how to introduce  what we would call a Poisson superalgebra structure on a Grassmann algebra and applied  that to obtain a Lagrangian classical counterpart of a quantum particle of spin ${1}/{2}\,$; in the second  paper, he started from a general algebraic formalism linking matrix calculus with nilpotent variables and applied it to constructing a Feynman integral over histories for fermionic systems. With hindsight, we may observe that Martin in these two works published together introduced  the integral over anticommuting variables. Strangely, he   applied  the name `integral' only for the functional case treated in~\cite{martin:feynman}. For the finite-dimensional case, he spoke about   cosets  \emph{modulo} total differentials in~\cite{martin:classical} or an unnamed `operation $S_{\lambda}$' in~\cite{martin:feynman}, $\lambda$ being a Grassmann algebra generator or, more generally, an abstract nilpotent variable. Martin did not consider transformations of variables and the corresponding properties of the integral. It is amazing that the remarkable works~\cite{martin:classical, martin:feynman} were not continued  and remained completely unnoticed. (Berezin learned about them only around 1976, see~\cite{berezin:marinov1}. He gives a very generous reference to~\cite{martin:classical} in the first sentence of~\cite{berezin:basis}.)

}

\smallskip
\subsection{From Grassmann algebras to supermanifolds.}
Berezin's calculus on a Grassmann algebra as constructed in~\cite{Ber-second1, BMV, Ber-second2, Ber-second3}  was not yet   supermathematics in the proper sense. Ordinary  variables and Grassmann algebra generators  were considered   in parallel but  still separately. There was no mixture of them nor transformations of Grassmann variables other than linear.
However, as Berezin described it later, \emph{``the striking coincidence of the main formulas of the operator calculus in the Fermi and Bose variants of the second quantization method}\dots      \emph{led to the idea of the possibility of  a generalization of all the main notions of analysis so that   generators of a Grassmann algebra would be on an equal footing with real or complex variables''}~\cite{berezin:basis, Ber-super-book}.

This was the \textbf{program of supermathematics}\,\footnote{Of course, the name came after the program was actually fulfilled.}.

The main steps of its implementation were as follows.

In~\cite{berezin:autgrass},  Berezin considered non-linear transformations of anticommuting variables in a clear departure from the standard viewpoint on the exterior algebra  as a $\ZZ$-graded object associated with  a linear space. Now  the emphasis is shifted to the algebra itself and the  transformations are supposed to preserve  only $\Z$-grading and  not necessarily  $\ZZ$-grading. Berezin studied the effect of such transformations on the integral over anticommuting variables and proved that there appears the \emph{inverse} of the determinant of the Jacobi matrix. This was a generalization of  the formula in~\cite{Ber-second} and  a  step towards the discovery of Berezinian. No mixture   with ``ordinary'' variables yet, but the whole logic leads in this direction.

Algebras generated by even and odd variables  appeared  in a joint paper of Berezin and G.~I.~Kac~\cite{berezin:andkac}, who introduced---in  modern language---formal Lie supergroups and Lie superalgebras and established  their 1-1 correspondence. They used the results of Milnor and Moore~\cite{milnor:andmoore}. It should be said that a version of Lie superalgebras where the $\Z$-grading arises as the reduction of a $\ZZ$-grading \emph{modulo} $2$ had been long familiar to topologists and differential geometers under the confusing name of ``graded Lie algebras''\,\footnote{A possible confusion with  the ordinary Lie algebras possessing  a grading.}. The understanding that graded (co)commutative Hopf algebras played the role of the corresponding group objects was topologists' folklore\,\footnote{As S.~P.~Novikov told to the writer of these words many years ago,  much of Milnor and Moore's paper had been part of folklore before its publication.}. In algebraic topology, Hopf algebras arise as homology or cohomology of topological spaces, so in that context $\ZZ$-grading is natural. Unlike that, the algebras  considered in~\cite{berezin:andkac} were supposed to play the role of  algebras of functions  and  the natural grading is  $\Z$.  Though~\cite{berezin:andkac} was devoted to the analogs of formal groups, the authors explained what the analog of a non-formal Lie group should be and gave  two examples: in modern language, the general linear supergroup $\GL(n|m)$ and the  diffeomorphism supergroup   $\Diff(\R{0|m})$.

\smallskip
{\small

(The notion of a Hopf algebra was discovered by   Milnor, motivated by the study of cohomology operations.   G.~I.~Kac, who should not be confused with V.~G.~Kac of Kac--Moody algebras,  independently came to a close notion, which he called a `ring group', working in representation theory. It  was instrumental for his generalization of the Pontrjagin duality and the Tannaka--Krein duality. Works of G.~I.~Kac, who died untimely in 1978, anticipated the discovery of quantum groups; incidentally, before quantum groups no good examples of Hopf algebras  that are neither commutative nor cocommutative were known. So the collaboration of Berezin and G.~I.~Kac on~\cite{berezin:andkac}  was not accidental.   Before~\cite{berezin:andkac}, Nijenhuis came very close to the concept of a Lie supergroup in deformation theory. Nijenhuis  used pairs consisting of a Lie superalgebra---of course, $\Z$ was $\ZZ$ {\emph{modulo}} $2$---and a Lie group corresponding to its even part, which is an ordinary Lie algebra. Such pairs are equivalent to Lie supergroups and are nowadays sometimes referred to by the name `Harish-Chandra pairs', borrowed from representation theory.)

}

\smallskip
After~\cite{berezin:autgrass, berezin:andkac}, everything seems  ready for the introduction of ``spaces'' for which the elements of $\Z$-graded algebras would be ``functions''.  But in fact it required a few  more years and some extra steps.

Such ``spaces'' remain implicit in paper~\cite{g.i.kac:andkoronkevich} by G.~I.~Kac and A.~I.~Koronkevich, submitted shortly after~\cite{berezin:andkac}, where   a superanalog of Frobenius theorem in the language of differential forms was stated and proved.

A preliminary step  was made in the   setting  of algebraic geometry. In paper~\cite{leites:spectra}, submitted in February 1973, Berezin's student D.~A.~Leites introduced a generalization of affine schemes (over a field) to the case of $\Z$-graded algebras. In particular, he introduced affine group schemes in this context and  defined their Lie algebras (in the sense of~\cite{berezin:andkac}, i.e, Lie superalgebras).

The missing ingredient---before differentiable supermanifolds  would  become pos\-sible---was transformations of variables mixing the ordinary variables with Grassmann generators. Berezin came to the idea of such transformations studying his integral: about the same time as~\cite{berezin:andkac} was written, he arrived at a formula for a general change of variables in the integral over a collection of anticommuting and ordinary variables.  According to Minlos~\cite{minlos}, the conjectural statement originally appeared  in 1971 in a letter to G.~I.~Kac. It contained, in particular, the notion of a `superdeterminant' (this name emerged only later\,\footnote{Now the term `Berezinian' and the notation $\Ber$ are universally adopted.}):
\begin{equation}\label{eq.defber}
    \sdet\!\begin{pmatrix}
            A & B \\
            C & D \\
          \end{pmatrix}:=  {\det(A-BD^{-1}C)}\,({\det D})^{-1}\,.
\end{equation}
Here the entries of $A$, $D$ are even and the entries of $B$, $C$, odd.
The change of variables formula reads (in the notation close to Berezin's own notation):
\begin{align}\label{eq.change}
    \int f(y,\h)\,d\h dy &=\int f(y(x,\x),\h(x,\x))\,J(x,\x)\,d\x dx\,,
\intertext{where}
\label{eq.jack}
    J(x,\x)&=\sdet\! \begin{pmatrix}
            \der{y}{x} & \der{\h}{x} \\[1mm]
            \der{y}{\x} & \der{\h}{\x} \\
          \end{pmatrix}\,.
\end{align}
The integral is extended to all values of the variables $y$ and the function $f(y,\h)$ must be vanishing at the infinity of $y$.

It is a curious fact that Berezin did not publish the definition of superdeterminant and the change of variables formula himself.  Berezin suggested a proof of~\eqref{eq.change} as a problem to his student V.~F.~Pakhomov, in whose paper~\cite{pakhomov},  submitted in December 1973, the above  formulas  first appeared in print\,\footnote{Without any particular name and notation for the function $\sdet$.}.

Allowing changes of variables in the Berezin  integral implied considering ordinary coordinates and Grassmann generators on an equal footing  as  generators of the algebra $\fun(\R{n})\otimes \Lambda(\R{m})$,  denoted $\mathfrak{B}_{n,m}$ in~\cite{pakhomov} (in modern language it is $\fun(\R{n|m})$). This was probably the final step towards supermanifolds.

An algebraic proof of the multiplicativity of ``Berezin's function''~\eqref{eq.defber} was given by  D.~A.~Leites in a short note ~\cite{leites:sdet}, submitted  in May 1974.

\subsection{Emergence of supersymmetric models and the explicit introduction of supermanifolds.}\label{subsec.supersym}
The   analysis in the previous subsection amply demonstrates that  Berezin's program had been mainly fulfilled by himself and his collaborators by around 1973. The notion of a supermanifold was  for them ``in the air'', though it had not appeared in the publications explicitly. The same can be said about supergroups. Still, according to Leites,  Berezin felt reluctant to publish the definition of a supermanifold and was forced to do so only in order not to lose the priority. So what happened?

The momentum came again from   theoretical physics and it was  \textbf{supersymmetry}.  Now this name is used very widely and sometimes outside of its precise original meaning, which is transformations of fields mixing the fermionic fields (usually describing `matter') with bosonic fields (usually describing `interaction').

In parallel with Berezin's work, the breakthrough  was preparing  in 1971--1974.
Supersymmetry  appeared, in the context of `dual-resonance models' (later, string theory), in Ramond~\cite{ramond:dual71} and Neveu--Schwarz~\cite{neveu:schwarz1971}; and, in the context of four-dimensional gauge theory, in Golfand--Likhtman~\cite{golfand:andlikhtman71}, Volkov \emph{et al.}~\cite{volkov:akulov72, volkov:akulov74, volkov:soroka74} (who explicitly quoted Berezin and G. I. Kac~\cite{berezin:andkac}), and finally in Wess--Zumino~\cite{wesszumino:supergauge74, wesszumino:lagrmodel74}, whose work resulted in an explosion. Physicists started to look around for mathematical foundations of the new theory.
Salam and Strathdee~\cite{salamstrathdee:supergauge} were the first to formulate the concept of a  `superspace' on an operational level.

In such a context, Berezin was forced to act quickly. Berezin and Leites published~\cite{berezin:andleites}. This paper contained the definition of a supermanifold as a local ringed space modeled on open domains of $\R{p}$ endowed with the $\Z$-graded algebra $\fun(\R{p})\otimes \Lambda(\R{q})$; the notions of morphisms of supermanifolds, subsupermanifolds and the direct products; the coordinate description by local charts and coordinate transformations; the notion of what we now call Berezin volume density and the construction of the Berezin integral over a supermanifold;  specifying supermanifolds by equations in $\R{n|m}$ and a conjecture that this may be possible for an arbitrary supermanifold (analog of Whitney's theorem); Lie supergroups (global) and their Lie superalgebras (so renamed from the `Lie algebras' of~\cite{berezin:andkac}). Quite a lot! Of course, the big work remained to elaborate the details and to make them available to the public.

It is time to say something about the terminology which involves the prefix `super-': supermanifold, superalgebra, supergroup, superspace...
In physical context, the term ``supersymmetry'' has a direct meaning as a ``superior symmetry'' exceeding other symmetries that keep bosons and fermions separate. In mathematics, the prefix `super-' should be understood as an abbreviation from supersymmetry,---as having something to do with physical supersymmetries.  Berezin himself did not overuse this prefix. It was done by others, and this made an unfortunate aftertaste. Nevertheless, such is the universally adopted terminology and there is no other choice but follow it.

\subsection{Berezin's work  on supermathematics in 1975--1980.}\label{subsec.berezinlatest}

There were several directions of Berezin's work after the introduction of supermanifolds.

The physical papers by Berezin and Marinov~\cite{berezin:marinov1, berezin:marinov2} were devoted to the description of spin by means of supermathematics. These works are still on the borderline with the previous period: they do not use the word `supermanifold' yet; the earlier paper~\cite{berezin:marinov1}  was submitted for publication just one month after~\cite{berezin:andleites}. There is an interesting historical material (mainly physical) in~\cite{berezin:marinov2}, in particular,  references to Martin~\cite{martin:classical, martin:feynman}.
Berezin came to this subject again  in a joint paper with V.~L.~Golo~\cite{berezin:golo} of 1980. It appeared only a few days before Berezin died; he could not see it published. One can also mention here  the posthumous publication~\cite{berezin:sigma} on a chiral supersymmetric sigma-model.

A central topic of Berezin's research in 1975--1980 was the theory of Lie superalgebras and Lie supergroups, especially their representations and invariants. Berezin's methods were global, geometric and analytic (e.g., used tools such as invariant integral) rather than infinitesimal.

In a short article~\cite{Ber-superupq}, Berezin studied the Lie supergroup $\U(p\,|q)$ and its unitary representations\,\footnote{Berezin's own notation for this Lie supergroup was $U(p,q)$ and this should not lead to a confusion with the ordinary Lie group of pseudounitary matrices.}. In particular, Berezin found the invariant integral and the radial parts of Laplace operators; he introduced   ``non-degenerate'' or ``typical'' irreducible representations and found their characters.

The method sketched in~\cite{Ber-superupq} was    elaborated and generalized in a series of
five preprints~\cite{berezin:itepla, berezin:iteplg, berezin:iteplapcas, berezin:iteprad, berezin:itepreps} of 1977. They contain   very interesting material; it would be fair to say that much of it has   unfortunately remained not well understood yet. These preprints,   originally published in a   small number of copies, were later included in the expanded English edition\,\footnote{We should be thankful to the editors for that. The quality of the English translation is sometimes poor: e.g., ``resultant'' is confused with ``result'', but it is still readable.} of~\cite{Ber-super-book}.

Two joint papers of Berezin and V.~S.~Retakh~\cite{berezin:superssla1, berezin:superssla2} of 1978 were devoted to classification of   Lie superalgebras whose even part is semisimple. (The classification of simple Lie superalgebras over $\CC$ was obtained by V.~G.~Kac around 1975, who interacted actively with Berezin at that time. Kac's classification remarkably brings forward superanalogs of classical matrix Lie algebras and the Lie algebras of vector fields, which is one more evidence for the ``naturalness'' of supergeometry.)

Berezin's last publication on representations of Lie supergroups was the paper with V.~N.~Tolstoy~\cite{berezin:itolstoy} dealing with a certain real compact form of the Lie supergroup $\OSp(1|2)$. (It appeared already after his death.)

Besides the Lie supergroups, Berezin actively worked on the general theory of supermanifolds. We should mention the expository preprint~\cite{berezin:itepsuperman} and the survey paper~\cite{berezin:basis} (both of 1979), and of course Berezin's work on a book on supermanifolds, which was incomplete at the time of his death. It was to appear only as a posthumous publication~\cite{Ber-super-book}, compiled and edited by his friends such as A.~A.~Kirillov and V.~P.~Palamodov.  (Palamodov, in particular, included there his own new results on the structure of supermanifolds.) The Russian version appeared  in 1983  and   the expanded English translation in 1987.

Three mathematical questions that attracted Berezin's attention are worth mentioning.

The first was the question about ``points of supermanifolds''. No doubt, the fact that ``functions'' on supermanifolds contain nilpotents makes it harder to understand them as compared to ordinary manifolds.  Supermanifolds cannot be treated as sets with some structure. For example, the supermanifold $\R{0|m}$, whose ``algebra of functions'' is  the Grassmann algebra with $m$ generators,  set-theoretically consists of a single point; clearly, the structure of $\R{0|m}$ cannot be attributed to this one-point set. At the same time, physicists working with ``superspace'' freely used ``points'' such as $(x^a,\xi^{\mu})$ with the odd coordinates $\xi^{\mu}$, whatever that could mean. Berezin's solution to that was in the introduction of an auxiliary Grassmann algebra $\mathfrak{G}(N)$ with a large or infinite number of generators $N$ and in considering, for each supermanifold $M$, the ordinary manifold $M(N)$ (of large dimension) obtained by replacing abstract even and odd coordinates in the coordinate transformations for $M$ by  elements of $\mathfrak{G}(N)$ of corresponding parities. The manifold $M(N)$, by construction, has a special structure called by Berezin a \emph{Grassmann-analytic} structure\,\footnote{This exactly means that the coordinate transformations on $M(N)$ are given by transformations of elements of the algebra $\mathfrak{G}(N)$ taken as   whole quantities---like complex numbers instead of their real and imaginary parts taken separately. This is a particular case of what is called a `manifold over an algebra'.}. If $N$ is large enough, the supermanifold $M$ can be recovered from $M(N)$ taken with this structure. At the same time, the `Grassmann-analytic manifold' $M(N)$ is a set-theoretic object and can be described by its points. Berezin  used this idea widely, in particular for representations of Lie superalgebras and Lie supergroups, which he replaced by (in his terminology)   their `Grassmann envelopes'.

Berezin's idea about the manifolds $M(N)$ and `Grassmann-analytic manifolds' in general contained the roots of  several later developments (some probably independent of him). If one fixes a Grassmann algebra  $\mathfrak{G}(N)$ and considers manifolds over it endowed with some class of `Grassmann-analytic functions', there is a temptation to   forget about supermanifolds defined as ringed spaces altogether. Two versions of this idea were put forward, by B.~DeWitt and by Alice Rogers, but in spite of all intuitive attractiveness, it was later found that its consistent development takes one back to the sheaf-theoretic approach to supermanifolds. Another option is not to fix  $\mathfrak{G}(N)$, but consider the manifolds $M(N)$ as a functor of $\mathfrak{G}(N)$. Together they represent the original supermanifold $M$. This was suggested by A.~S.~Schwarz. It allows to consider objects more general than supermanifolds.

The second question concerned general classification of supermanifolds. It is obvious that the case in which the coordinate changes for a given supermanifold do not mix odd variables with the even variables and the odd variables transform linearly is the simplest and it corresponds to a vector bundle over an ordinary manifold. The question is how general this case is, i.e., whether it is always possible to reduce coordinate transformations to this simple form by a choice of atlas. The answer is, yes,--- for smooth supermanifolds. This statement is often referred to as the Batchelor theorem after Marjorie~Batchelor who proved it in~1979. No doubt that Berezin knew it independently: he mentions it in~\cite{berezin:itepsuperman, berezin:basis}. As for complex-analytic supermanifolds, the answer is, no; there are obstructions. The corresponding theory and examples of ``non-retractable'' complex-analytic supermanifolds are due to V.~P.~Palamodov  (a slightly different approach was developed by Yu.~I.~Manin).

Finally, the third question concerned integration on supermanifolds and differential forms. It is clear that Berezin's transformation law for the element of volume is different from what one gets for the differentials of coordinates defined in a straightforward way.  So  on supermanifolds, differential forms and integration theory seem to split. The problem was tackled by J.~N.~Bernstein and D.~A.~Leites, who introduced `integral forms'~\cite{berl:int}, incorporating volume elements, as a replacement of differential forms for the purpose of integration and   `pseudodifferential forms'~\cite{berl:pdf} as not necessarily polynomial expressions in differentials (which also opens way for integration).

Berezin in~\cite{berezin:difforms} introduced a further generalization of the Bernstein--Leites pseudodifferential forms, studied their  duality transformations and sketched a Weil-type construction of characteristic classes. Berezin's aim was future application to super gauge theory. Paper~\cite{berezin:difforms} seems to be Berezin's last paper on supermanifolds. It is worth mentioning that~\cite{berezin:difforms} contained a construction very close to what is now known as the `homological interpretation' of Berezinian.

\subsection{Influence. Later developments.}

The influence of Berezin's work on supermathematics remains different in physics and in pure mathematics. Physicists have completely absorbed the idea of working with supermanifolds. For them, supergeometry is a tool on the same footing as tensor calculus: physicists use it without even noticing it. Unlike that, in pure mathematics, Berezin's ideas have spread far less widely. Supermanifolds for many remain something exotic (except for those directly working in supergeometry). Quite characteristic is that representations of  Lie superalgebras became a well-established area, but those working in it rarely consider Lie supergroups or turn to global methods used by Berezin. No doubt, the landscape would be quite different, had Berezin not died in 1980.    However, the situation is slowly changing. ``Supermethods'' start to spread in differential geometry. Of course, this development is more significant in areas closer to or more influenced by physics. Two Fields medals awarded in 1990 and 1998, to E.~Witten and M.~L.~Kontsevich, respectively, were related with works where supergeometry played a role. (Morse theory and differential forms, in the case of Witten, and deformation theory and quantization of Poisson manifolds, in the case of Kontsevich.)

One general trend worth mentioning is a certain shift from ``supersymmetry'' (roughly, transformations squaring to ordinary symmetries) to ``BRST-symmetry'' and ``$Q$-manifolds'' (where, roughly, there are transformations with square zero). A central role has been played here by the Batalin--Vilkovisky formalism in quantum field theory~\cite{bv:perv, bv:vtor, bv:closure} and its modern generalizations. Geometrically, that means considering supermanifolds endowed with an odd symplectic structure and odd Laplacians on them. (The study of such geometry was pioneered by H.~M.~Khudaverdian~\cite{hov:deltabest}, see also~\cite{tv:laplace1, tv:laplace2}.)

Another trend is the growth of importance of \emph{graded manifolds}\,---\,not in the sense synonymous with supermanifolds as the usage in the early period  sometimes was\,\footnote{In the Western literature, when the foundations of supermanifolds were thought to be not fully established yet, `supermanifolds' were often used for DeWitt or Rogers's versions of manifolds over Grassmann algebras, while `supermanifolds' in Berezin's sense were called `graded manifolds'.},\,---\, but  meaning supermanifolds endowed with extra $\ZZ$- or $\ZZ_+$-grading, which in physics may be for example, `ghost number'. If one recalls topologists' $\ZZ$-graded algebras and the replacement
of $\ZZ$-grading by $\Z$-grading as a step in development of supermathematics, as
described above, the reintroduction of $\ZZ$-gradings (but now as additional structure)
completes the circle, but at a higher level.

%\smallskip
We would like to finish this section by two interesting pieces of mathematics related with Berezin integration and Berezinian (superdeterminant), which were discovered after Berezin.

In the previous subsection, we considered the works   on integration theory and differential forms by Bernstein--Leites and by Berezin~\cite{berezin:difforms}. In  search of objects suitable for integration over (multidimensional) paths or surfaces in supermanifolds, a variational approach to ``forms on supermanifolds'' was developed (Th.~Th.~Voronov and A.~V.~Zorich, see~\cite{tv:git}; building on earlier works by A.~S.~Schwarz and his students): analogs of forms were constructed as Lagrangians satisfying certain restrictions. An amazing fact discovered along these lines and not fully understood yet is the following link with integral geometry in the sense of Gelfand--Gindikin--Graev: the equation of the form\,\footnote{\,The tilde over an index denotes the parity of the corresponding  variable.}
\begin{equation}\label{eq.ours}
    \dder{f}{w_a^i}{w_b^j}+(-1)^{\itt\jtt +\at(\itt+\jtt)}\dder{f}{w_a^j}{w_b^i}=0\,,
\end{equation}
for a function of a rectangular supermatrix $\|w_a^i\|$, arises in the de Rham theory on supermanifolds as a condition replacing skew-symmetry and multilinearity (see~\cite{tv:git}) and at the same time it is a generalization of `hypergeometric equations' in the sense of Gelfand (the odd-odd part of~\eqref{eq.ours} is the F.~John equation arising in relation with the Radon transform).

Another beautiful development related with the notion of Berezinian is as follows. Th.~Schmitt~\cite{schmitt:ident} discovered that the expansion of Berezinian leads to exterior powers:
\begin{equation*}
    \Ber(1+ zA)=1+z\str A +z^2\str \Lambda^2A +\ldots\,.
\end{equation*}
Here $A$ is an even supermatrix,
\begin{equation*}
    A=\begin{pmatrix}
        A_{00} & A_{01} \\
        A_{10} & A_{11} \\
      \end{pmatrix},
\end{equation*}
$\Lambda^kA$ stands for its action on the $k$th exterior power and
$\str$ denotes the supertrace: $\str A=\tr A_{00}-\tr A_{11}$. As it was
found in~\cite{tv:ber}, by comparing the expansions of $\Ber(1+ zA)$
at zero and at infinity one arrives at certain universal recurrence
relations satisfied by the differences of the respective
coefficients. In particular, for a $p\,|q$ $\times$ $p\,|q$ matrix,
there are   relations
\begin{equation*}
\begin{vmatrix}
      c_k(A) & \dots & c_{k+q}(A) \\
      \dots & \dots & \dots \\
      c_{k+q}(A) & \dots & c_{k+2q}(A) \\
    \end{vmatrix}=0\,,
\end{equation*}
where $c_k(A)=\str \Lambda^kA$, satisfied for all $k> p-q$. (This
replaces the vanishing of the $k$th exterior powers for an
$n$-dimensional space with $k>n$). Similar relations hold in the
Grothendieck ring of a general linear supergroup, and there is a
formula
\begin{equation*}
    \Ber A%=\frac{\begin{vmatrix}
    %c_{p-q} & \ldots & c_p \\
     % \ldots & \ldots & \ldots \\
    %  c_p & \ldots & c_{p+q} \\
    %\end{vmatrix}}{\begin{vmatrix}
    %  c_{p-q+2} & \ldots & c_{p+1} \\
     % \ldots & \ldots & \ldots \\
   %   c_{p+1} & \ldots & c_{p+q} \\
    %\end{vmatrix}}%
    =\frac{|\,c_{p-q}(A)\ldots c_p(A)\,|_{q+1}}{|\,c_{p-q+2}(A)\ldots c_{p+1}(A)\,|_{q}}\,,
\end{equation*}
with Hankel's determinants at the top and at the bottom, expressing
Berezinian as the ratio of polynomial invariants\,\footnote{\,In the
  definition, $\Ber A=\det(A_{00}-A_{01}A_{11}^{-1}A_{10})\det
  A_{11}^{-1}$, neither the numerator nor the denominator of the
  fraction are invariant.}.

\bigskip
\subsection*{Acknowledgment}
We thank S.~G.~Gindikin for discussions.   A.~V.~Karabegov
and Yu.~A.~Neretin acknowledge the support by the NSF and RFBR
respectively that made their participation  in the XXXth
Bia{\l}owie\.{z}a conference possible. (The NSF  award number: 1124929.) Neretin's work was also partially supported by the FWF grant, Project 22122.

%\bibliographystyle{hplain} \bibliography{geometry} \end{document}

\begin{thebibliography}{200}


\bibitem{JEA} J.~E.~Andersen and J.~Blaavand. Asymptotics of Toeplitz
  operators and applications in TQFT. {\it Travaux Mathematiques} {\bf
    19} (2011), 167--201.

\bibitem{Bar} V.~Bargmann. On a Hilbert space of analytic functions
  and an associated integral transform.  \emph{Comm. Pure Appl. Math.}
  \textbf{14} (1961), 187--214.


\bibitem{bv:perv} I.~A. Batalin and G.~A. Vilkovisky.  \newblock Gauge
  algebra and quantization.  \newblock {\em Phys. Lett.} \textbf{102B}
  (1981), 27--31.

\bibitem{bv:vtor} I.~A. Batalin and G.~A. Vilkovisky.  \newblock
  Quantization of gauge theories with linearly dependent generators.
  \newblock {\em Phys. Rev.} \textbf{D28} (1983), 2567--2582.

\bibitem{bv:closure} I.~A. Batalin and G.~A. Vilkovisky.  \newblock
  Closure of the gauge algebra, generalized {Lie} equations and
  {Feynman} rules.  \newblock {\em Nucl. Phys. B} \textbf{234} (1984),
  106--124.


\bibitem{BFFLS} F.~Bayen, M.~Flato, C.~Fronsdal, A.~Lichnerowicz, and
  D.~Sternheimer. Deformation theory and quantization. I. Deformations
  of symplectic structures. {\it Ann. Physics} {\bf 111} (1978),
  no. 1, 61--110.


\bibitem{L1} F.~A.~Berezin.  Laplace operators on semisimple Lie
  groups. \emph{Dokl. Akad. Nauk SSSR (N.S.)} \textbf{107} (1956),
  9--12.

\bibitem{L2} F.~A.~Berezin.  Representation of complex semisimple Lie
  groups in Banach space.  \emph{Dokl. Akad. Nauk SSSR (N.S.)}
  \textbf{110} (1956), 897--900.


\bibitem{L3} F.~A.~Berezin.  Laplace operators on semisimple Lie
  groups and on certain symmetric spaces. (Russian) \emph{Uspehi
    Mat. Nauk (N.S.)} \textbf{12} (1957), no. 1(73),
  152--156. (Translated in: Amer. Math. Soc. Transl. (2) \textbf{16}
  (1960), 364--369.)

\bibitem{Ber-main} F.~A.~Berezin.  Laplace operators on semisimple Lie
  groups. \emph{Trudy Moskov. Math. Obshchestva.} \textbf{6} (1957),
  371--463.  (Translated in: Amer. Math. Soc. Transl. (2) \textbf{21}
  (1962), 239--339.)


\bibitem{Ber-second1} F.~A.~Berezin.  Canonical operator
  transformation in representation of secondary quantization.
  \emph{Dokl. Akad. Nauk SSSR} \textbf{137}, 311--314 (Russian);
  translated as \emph{Soviet Physics Dokl.} \textbf{6} (1961),
  212--215.

\bibitem{Ber-second2} F.~A.~Berezin.  Canonical transformations in the
  second quantization representation. (Russian) \emph{Dokl. Akad. Nauk
    SSSR} \textbf{150} (1963), 959--962.

\bibitem{Ber-second3} F.~A.~Berezin.  Operators in the representation
  of secondary quantization. \emph{Dokl. Akad. Nauk SSSR}
  \textbf{154}, 1063--1065 (Russian); translated as \emph{Soviet
    Physics Dokl.} \textbf{9} (1964), 142--144.


\bibitem{Ber-second} F.~A.~Berezin. \emph{The Method of Second
    Quantization.} Nauka, Moscow, 1965. Tranlation: Academic Press,
  New York, 1966. (Second edition, expanded: M.~K.~Polivanov, ed.,
  Nauka, Moscow, 1986.)

\bibitem{Ber-some-remarks} F.~A.~Berezin. Some remarks on the
  representations of commutation relations. \emph{Uspehi Mat. Nauk}
  \textbf{24}, no. 4 (148) (1969), 65--88; Russian
  Math. Surv. \textbf{24}, No.4, 65--88.

\bibitem{Ber-letter} F.~A.~Berezin.  Letter to the editor.
  \emph{Trudy Moskov. Mat. Obshch.} \textbf{12} (1963), 453--466.
% Est' anglijskij perevod.

\bibitem{berezin:autgrass} F.~A. Berezin.  \newblock Automorphisms of
  a {G}rassmann algebra.  \newblock {\em Mathematical Notes},
  \textbf{1} (1967), 180--184.


\bibitem{Ber5} F.~A. Berezin. Some remarks about the associative
  envelope of a Lie algebra. {\it Funct. Anal. Appl.} {\bf 1} (1967),
  91--102.

\bibitem{Ber6} F.~A. Berezin. Non-Wiener path integrals. (Russian)
  \emph{Teoret. Mat. Fiz.} {\bf 6} (1971), no. 2, 194--212.

\bibitem{Ber7} F.~A. Berezin. Wick and anti-Wick operator
  symbols. {\it Mat. Sb. (N.S.)} {\bf 15} (1971), 577--606.

\bibitem{Ber4} F.~A. Berezin. Covariant and contravariant symbols of
  operators. {\it Math. USSR-Izv.} {\bf 36} (1972), 1134--1167.
\bibitem{Ber9} F.~A. Berezin. Convex operator functions. {\it Mat. Sb. (N.S.)} {\bf 17} (1972), 269--278.

\bibitem{Ber-hol1} F.~A.~Berezin.  Quantization in complex bounded
  domains. (Russian) \emph{Dokl. Akad. Nauk SSSR} \textbf{211} (1973),
  1263--1266.

\bibitem{Ber1} F.~A. Berezin. Quantization. \emph{Math. USSR-Izv.}
  {\bf 38} (1974), 1116--1175.

\bibitem{Ber3} F.~A. Berezin. General concept of quantization. {\it
    Comm. Math. Phys.} {\bf 40} (1975), 153--174.


\bibitem{Ber2} F.~A. Berezin. Quantization in complex symmetric
  spaces. {\it Math. USSR-Izv.} {\bf 39} (1975), 363--402.
%\bibitem{Ber-hol2} F.~A.~Berezin.  Quantization in complex symmetric spaces. \emph{Izv. Akad. Nauk SSSR Ser. Mat.} \textbf{39} (1975), no. 2, 363�-402.

\bibitem{Ber-superupq} F.~A.~Berezin.  Representations of the
  supergroup $U(p,q)$.  \emph{Funct. Anal.  Appl} \textbf{10}(3)
  (1976), 221--223.  (Transl. from:
%\newblock Representations of the supergroup {$U(p,q)$}.
\newblock {\em Funkcion. Anal. i Prilo\v zen.}, \textbf{10}(3) (1976), 70--71.)


\bibitem{berezin:itepla} F.~A. Berezin.  \newblock Lie superalgebras.
  \newblock Preprint ITEP-66, 1977.

\bibitem{berezin:iteplg} F.~A. Berezin.  \newblock Lie supergroups.
  \newblock Preprint ITEP-78, 1977.

\bibitem{berezin:iteplapcas} F.~A. Berezin.  \newblock Laplace-Casimir
  operators (general theory).  \newblock Preprint ITEP-77, 1977.

\bibitem{berezin:iteprad} F.~A. Berezin.  \newblock The radial parts
  of the Laplace-Casimir operators on Lie supergroups $U(p,q)$ and
  $C(m,n)$.  \newblock Preprint ITEP-75, 1977.


\bibitem{berezin:itepreps} F.~A. Berezin.  \newblock Construction of
  representations of Lie supergroups $U(p,q)$ and $C(m,n)$.  \newblock
  Preprint ITEP-76, 1977.


\bibitem{berezin:itepsuperman} F.~A. Berezin.  \newblock
  Supermanifolds.  \newblock Preprint ITEP, 1979.


\bibitem{Ber10} F.~A. Berezin. The connection between covariant and
  contravariant symbols of operators on classical complex symmetric
  spaces. (Russian) \emph{Dokl. Akad. Nauk SSSR} \textbf{241} (1978),
  no. 1, 15--17.
%\bibitem{Ber-hol3} F.~A.~Berezin.  The connection between covariant and contravariant symbols of operators on classical complex symmetric spaces. (Russian) \emph{Dokl. Akad. Nauk SSSR} \textbf{241} (1978), no. 1, 15-�17.


\bibitem{berezin:basis} F.~A. Berezin.  \newblock The mathematical
  basis of supersymmetric field theories.  \newblock {\em Soviet
    J. Nuclear Phys.}, \textbf{29}(6) (1979), 1670--1687.

\bibitem{berezin:difforms} F.~A. Berezin.  \newblock Differential
  forms on supermanifolds.  \newblock {\em Soviet J. Nuclear Phys.},
  \textbf{30}(4) (1979), 605--609.


\bibitem{Ber8} F.~A. Berezin. Feynman path integrals in a phase
  space. (Russian) \emph{Soviet Phys. Uspekhi} \textbf{132} (1980),
  no. 3, 497--548.


\bibitem{Ber-super-book} F.~A.~Berezin. \emph{Introduction to algebra
    and analysis with anticommuting variables.} V.~P.~Palamodov, ed.,
  Moscow State University Press, Moscow, 1983. Expanded transl. into
  English: {\it Introduction to superanalysis.} A.~A.~Kirillov, ed.,
  D. Reidel, Dordrecht, 1987.


\bibitem{BG}%% formatted in the required style
  F.~A.~Berezin and I.~M.~Gelfand. Some remarks on the theory of
  spherical functions on symmetric Riemannian manifolds. (Russian)
  \textit{Trudy Moskov. Mat. Obshch.} \textbf{5} (1956), 311--351.


\bibitem{bggn} F.~A.~Berezin, I.~M.~Gelfand, M.~I.~Graev. Group
  representations.  \emph{Uspehi Mat. Nauk (N.S.)} \textbf{11} (1956),
  no. 6(72), 13--40.  (Amer. Math. Soc. Transl. (2) \textbf{16}
  (1960), 325--353.)

\bibitem{berezin:golo} F.~A. Berezin and V.~L. Golo.  \newblock
  Supersymmetric model of several classical particles with
  spin. \newblock (Russian) \newblock {\em JETP Letters}
  \textbf{32}(1) (1980), 82--84.


\bibitem{berezin:sigma} F.~A. Berezin, V.~L. Golo, and B.~A. Putko.
  \newblock Reduction of a chiral supersymmetric {$\sigma $}-model.
  \newblock {\em Vestnik Moskov. Univ. Ser. I Mat. Mekh.}, (3) (1982),
  16--19, 109.

\bibitem{berezin:andkac} F.~A. Berezin and G.~I. Kac.  \newblock Lie
  groups with commuting and anticommuting parameters.  \newblock {\em
    Mathematics of the {USSR-Sbornik}}, \textbf{11}(3) (1970),
  311--325.


\bibitem{BK} F.~A.~Berezin and F.~I.~Karpelevich. Zonal spherical
  functions and Laplace operators on some symmetric
  spaces. \emph{Dokl. Akad. Nauk SSSR} \textbf{118} (1958), 9--12.

\bibitem{berezin:andleites} F.~A. Berezin and D.~A. Leites.  \newblock
  Supermanifolds.  \newblock {\em Soviet Math. Dokl.} \textbf{16}(5)
  (1975), 1218--1222.

\bibitem{berezin:marinov1} F.~A. Berezin and M.~S. Marinov.  \newblock
  Classical spin and the Grassmann algebra.  \newblock {\em JETP
    Letters} \textbf{21}(11) (1975), 678--680.

\bibitem{berezin:marinov2} F.~A. Berezin and M.~S. Marinov.  \newblock
  Particle spin dynamics as the Grassmann variant of classical
  mechanics.  \newblock {\em Annals of Physics} \textbf{104} (1977),
  336--362.


\bibitem{BMV} F.~A.~Berezin, R.~A.~Minlos, and L.~D.~Faddeev. Some
  problems of quantum-mechanical systems with large number of the
  degrees of freedom.
%��������� �������������� ������� ��������� �������� ������ �
%������� ������ �������� ������� // ����� 4-�� ��������. ��-���. ������. �., 1964. �. 2. �. 532�541.
  \emph{Proceedings of 4-th All-Union Congress of Mathematics (1961)},
  Moscow, 1964, vol. 2, 532--541.

\bibitem{berezin:kdvfa} F.~A. Berezin and A.~M. Perelomov.  \newblock
  Group-theoretic interpretation of equations of the
  {K}orteweg-de\thinspace {V}ries type.  \newblock {\em
    Funktsional. Anal. i Prilozhen.} \textbf{14}(2) (1980), 50--51.

\bibitem{berezin:kdvcmp} F.~A. Berezin and A.~M. Perelomov.  \newblock
  Group-theoretical interpretation of the
  {K}orteweg\thinspace-\thinspace de {V}ries type equations.
  \newblock {\em Comm. Math. Phys.} \textbf{74}(2) (1980), 129--140.

\bibitem{berezin:superssla1} F.~A. Berezin and V.~S. Retah.  \newblock
  The structure of {L}ie superalgebras with a semisimple even part.
  \newblock {\em Funkcional. Anal. i Prilo\v zen.} \textbf{12}(1)
  (1978), 64--65.

\bibitem{berezin:superssla2} F.~A. Berezin and V.~S. Retah.  \newblock
  Construction of {L}ie superalgebras with a semisimple even part.
  \newblock {\em Vestnik Moskov. Univ. Ser. I Mat. Mekh.} (5) (1978),
  63--67.

\bibitem{berezin:charcl} F.~A. Berezin and V.~S. Retakh.  \newblock A
  method of computing characteristic classes of vector bundles.
  \newblock {\em Rep. Math. Phys.}  \textbf{18}(3) (1980), 363--378.


\bibitem{berezin:ishubin1} F.~A. Berezin and M.~A. Shubin.  \newblock
  Symbols of operators and quantization.  \newblock In {\em Hilbert
    space operators and operator algebras ({P}roc.
    {I}nternat. {C}onf., {T}ihany, 1970)}, pages
  21--52. Colloq. Math. Soc.  J\'anos Bolyai, No. 5. North-Holland,
  Amsterdam, 1972.

\bibitem{berezin:ishubin} F.~A. Berezin and M.~A. Shubin.  \newblock
  {\em The {S}chr\"odinger equation}.  \newblock Moscow State
  University Press, Moscow, 1983. Transl. into English: \newblock
  Kluwer, Dordrecht, 1991.

\bibitem{berezin:itolstoy} F.~A. Berezin and V.~N. Tolstoy.  \newblock
  The group with {G}rassmann structure {${\rm UOSP}(1.2)$}.  \newblock
  {\em Comm. Math. Phys.}  \textbf{78}(3) (1980/81), 409--428.


\bibitem{berl:int} J.~N. Bernstein and D.~A. Leites.  \newblock
  Integral forms and {Stokes} formula on supermanifolds.  \newblock
  {\em Funk. Anal. Pril.} \textbf{11}(1) (1977), 55--56.

\bibitem{berl:pdf} J.~N. Bernstein and D.~A. Leites.  \newblock How to
  integrate differential forms on supermanifolds.  \newblock {\em
    Funk. Anal. Pril.} \textbf{11}(3) (1977), 70--71.

\bibitem{berezin:obituary-umn} N.~N. Bogolyubov, I.~M. Gelfand,
  R.~L. Dobrushin, A.~A. Kirillov, M.~G. Krein, D.~A. Leites,
  R.~A. Minlos, Ya.~G. Sinai, and M.~A. Shubin.  \newblock Feliks
  {A}leksandrovich {B}erezin. {O}bituary.  \newblock {\em Uspekhi
    Mat. Nauk} \textbf{36}(4(220)) (1981), 185--190.

\bibitem{BMS} M.~Bordemann, E.~Meinrenken, and
  M.~Schlichenmaier. Toeplitz quantization of K\"ahler manifolds and
  $gl(n), n \to\infty$, limits. {\it Comm. Math. Phys.} {\bf 165}
  (1995), 281--296.

\bibitem{BKLR} D. Borthwick, S. Klimek, A. Lesniewski, and M. Rinaldi. \newblock Super Toeplitz operators and non-perturbative quantization of supermanifolds, \newblock {\em Commun. Math. Phys.} \textbf {153} (1993), 49 -- 76.

\bibitem{BMG} L.~Boutet de Monvel and V.~Guillemin. \emph{The spectral
    theory of Toeplitz operators. } {Ann. Math. Stud.} {\bf 99},
  Princeton University Press, Princeton, 1981.

\bibitem{CGR1} M.~Cahen, S.~Gutt, and J.~Rawnsley. Quantization of
  K\"ahler manifolds I: Geometric interpretation of Berezin's
  quantization. {\it JGP} {\bf 7} (1990), 45--62.

\bibitem{CGR2} M.~Cahen, S.~Gutt, and J.~Rawnsley.  Quantization of
  K\"ahler manifolds II. {\it Trans. Amer. Math. Soc.} {\bf 337}
  (1993), 73--98.

\bibitem{LCh} L.~Charles. Berezin-Toeplitz operators, a semi-classical
  approach. \emph{Comm. Math. Phys.} {\bf 239} (2003), 1--28.

\bibitem{E} M.~Engli\v{s}.  Weighted Bergman kernels and
  quantization. {\it Comm. Math. Phys.} {\bf 227} (2002), 211--241.


\bibitem{Fel} J.~Feldman. Equivalence and perpendicularity of Gaussian
  processes. \emph{Pac. J. Math.} \textbf{8} (1958), 699--708.

\bibitem{Fri} K.~O.~Friedrichs. \emph{Mathematical aspects of the
    quantum theory of fields. } Interscience Publishers, Inc., New
  York, 1953.

\bibitem{Ful} W.~Fulton. Eigenvalues, invariant factors, highest
  weights, and Schubert calculus. \emph{Bull. Am. Math. Soc., New
    Ser.} \textbf{37}, No.3 (2000), 209--249.


\bibitem{GN} I.~M.~Gelfand and M.~A.~Naimark. {\it Unitary
    representations of classical groups.}  Trudy Mat.Inst. Steklov.,
  v.36 (1950) (Russian); German translation: I.~M.~Gelfand and
  M.~A.~Neumark.  {\it Unitare {D}arstellungen der klassischen
    {G}ruppen}, Akademie-{V}erlag, Berlin, 1957.

\bibitem{God} R.~Godement.  A theory of spherical
  functions. I. T\emph{rans. Am. Math. Soc.} \textbf{73} (1952),
  496--556.

\bibitem{golfand:andlikhtman71} Yu.~A.~Golfand and E.~P.~Likhtman.
  \newblock Extension of the algebra of {Poincar\'{e}} group
  generators and violation of {P} invariance.  \newblock {\em JETP
    Letters} \textbf{13} (1971), 323--326.

\bibitem{GAN} A. El Gradechi and L. M. Nieto. \newblock Supercoherent States, Super K\" ahler Geometry and Geometric Quantization. \newblock \emph{Commun.Math.Phys.} \textbf{175} (1996), 521--564.

\bibitem{Gui} A.~Guichardet and D.~Wigner. Sur la cohomologie r\'eelle
  des groupes de Lie simples r\'eels.  \emph{Ann. Sci. \'Ecole
    Norm. Sup.} (4) 11 (1978), no. 2, 277--292.

\bibitem{Harish1} Harish-Chandra. Representations of semisimple Lie
  groups, II.  \emph{Trans. Amer. Math. Soc.} \textbf{76} (1953),
  26--65.


%\bibitem{Harish1.5}
%Harish-Chandra
%{\it  Representations of semisimple Lie groups.  III.}
%  Trans. Am. Math. Soc. 76,
%              234-253 (1954).


\bibitem{Harish2} Harish-Chandra.  The characters of semisimple Lie
  groups.  \emph{Trans. Amer. Math. Soc.}  \textbf{83} (1956),
  98--163.


\bibitem{Heck} G.~Heckman and H.~Schlichtkrull.  \emph{Harmonic
    analysis and special functions on symmetric spaces.} (English)
  Perspectives in Mathematics. 16. Orlando, FL: Academic Press,
  Inc. xii, 225 p.

 \bibitem{Huck} A.~Huckleberry and M.~Kalus. Selected results on Lie
   supergroups and their radial operators.  Preprint:
   \texttt{arXiv:1012.5233}.

\bibitem{g.i.kac:andkoronkevich} G.~I. Kac and A.~I. Koronkevi{\v{c}}.
  \newblock The {F}robenius theorem for functions of commuting and
  anticommuting arguments.  \newblock {\em Funkcional. Anal. i Prilo\v
    zen.}  \textbf{5}(1) (1971), 78--80.

\bibitem{CMP1} A.~Karabegov. Deformation quantizations with separation
  of variables on a K\"ahler manifold. {\it Comm. Math. Phys.} {\bf
    180} (1996), 745--755.
%\bibitem{CMP2} Karabegov, A.: Pseudo-K\"ahler quantization on flag manifolds. {\it Commun. Math. Phys.} {\bf 200} (1999), 355--379.
\bibitem{KSch} A.~Karabegov and M.~Schlichenmaier. Identification of
  Berezin-Toeplitz deformation quantization. {\it J. reine
    angew. Math.} {\bf 540} (2001), 49--76.

\bibitem{hov:deltabest} H.~M. Khudaverdian.  \newblock Geometry of
  superspace with even and odd brackets.  \textit{J. Math. Phys.}
  \textbf{32} (1991), 1934--1937. (Preprint of the {Geneva
    University}, {UGVA-DPT} 1989/05-613, 1989.)

\bibitem{tv:laplace1}
H.~M.~Khudaverdian and Th. Th.~Voronov.
\newblock On odd {Laplace} operators.
\newblock {\em Lett. Math. Phys.} \textbf{62} (2002), 127--142.

\bibitem{tv:laplace2}
H.~M.~Khudaverdian and Th. Th.~Voronov.
\newblock On odd {Laplace} operators. {II}.
\newblock In book: \emph{Geometry, Topology and Mathematical Physics. S. P. Novikov�s
seminar: 2002 - 2003}, V.~M. Buchstaber and I.~M. Krichever, editors,  \emph{Amer. Math. Soc. Transl. (2)},
Vol. 212, 2004, pp.179--205.

\bibitem{tv:ber} H.~M. Khudaverdian and Th.~Th. Voronov.  \newblock
  Berezinians, exterior powers and recurrent sequences.  \newblock
  {\em Lett. Math. Phys.} \textbf{74}(2) (2005), 201--228.


\bibitem{Klya} A.~A.~Klyachko. Stable bundles, representation theory
  and Hermitian operators. \emph{Selecta Math., New Ser.} \textbf{4},
  No. 3 (1998), 419--445.

\bibitem{leites:spectra} D.~A. Leites.  \newblock Spectra of
  graded-commutative rings.  \newblock {\em Uspehi Mat. Nauk}
  \textbf{29}(3(177)) (1974), 209--210.


\bibitem{leites:sdet} D.~A. Leites.  \newblock A certain analogue
  of the determinant.  \newblock {\em Uspehi Mat. Nauk}
  \textbf{30}(3(183)) (1975), 156.

\bibitem{MM} X.~Ma and G.~Marinescu. \emph{Holomorphic Morse
    inequalities.} Progress in Mathematics, {\bf 254} (2007),
  Birkh\"auser Verlag, Basel.

\bibitem{martin:classical} J.~L. Martin.  \newblock Generalized
  classical dynamics, and the ``classical analogue'' of a {F}ermi
  oscillator.  \newblock {\em Proc. Roy. Soc. London. Ser. A}
  \textbf{251} (1959), 536--542.


\bibitem{martin:feynman} J.~L. Martin.  \newblock The {F}eynman
  principle for a {F}ermi system.  \newblock {\em
    Proc. Roy. Soc. London. Ser. A} \textbf{251} (1959), 543--549.

\bibitem{milnor:andmoore} J.~W. Milnor and J.~C. Moore.  \newblock On
  the structure of {H}opf algebras.  \newblock {\em Ann. of Math. (2)}
  \textbf{81} 1965, 211--264.


\bibitem{minlos} R.~A. Minlos.  \newblock Felix {A}lexandrovich
  {B}erezin (a brief scientific biography).  \newblock {\em
    Lett. Math. Phys.} \textbf{74}(1) (2005), 5--19.

\bibitem{Mol1} V.~F. Molchanov. \newblock Quantization on the imaginary Lobachevskii plane. \newblock{\em
Funkts. Anal. Prilozh.,} \textbf{14:2} (1980),  73–74.

\bibitem{Mol2} V.~F. Molchanov. \newblock Quantization on para-Hermitian symmetric spaces. \newblock{\em Amer.
Math. Soc. Transl.}, Ser. 2, \textbf{175} (1996), ({\em Adv. in Math. Sci.}, 31), 81--95.

\bibitem{Mor} C.~Moreno. $*$-Products on some K\"ahler manifolds. {\it
    Lett. Math. Phys.} {\bf 11} (1986), 361--372.

\bibitem{Ner-Ber} Yu.~A.~Neretin. Plancherel formula for Berezin
  deformation of $L^2$ on Riemannian symmetric
  space. \emph{J. Funct. Anal.} \textbf{189}, No. 2 (2002), 336--408.

\bibitem{Ner} Yu.~A.~Neretin. ``The method of second quatization'': a
  view 40 years after. \\
  In book: \emph{Recollections of Felix Alexandovich Berezin,
the founder of supermathematics},
Moscow, MCCME publishers, 2009 (Russian), available at
\url{http://www.mat.univie.ac.at/~neretin/zhelobenko/berezin.pdf}.
Translation into French by C.~Roger (with participation of O.~Kravchenko, D.~Millionschikov and A.~Kosyak) is available at
\url{http://hal.archives-ouvertes.fr/docs/00/47/84/76/PDF/neretin.pdf}.


\bibitem{Ner-gauss} Yu.~A.~Neretin.  {\it Lectures on Gaussian
    integral operators and classical groups.}  EMS Series of Lectures
  in Mathematics. Z\"urich: European Mathematical Society, 2011, xii,
  559 p.

\bibitem{neveu:schwarz1971} A.~Neveu and J.~H. Schwarz.  \newblock
  Factorizable dual model of pions.  \newblock {\em Nuclear Phys. B}
  \textbf{31} (1971), 86--112.

\bibitem{PO} M.~A.~Olshanetsky and A.~M.~Perelomov.  Completely
  integrable Hamiltonian systems connected with semisimple Lie
  algebras.  \emph{Inv. Math.}  \textbf{37} (1976), 93--108.

\bibitem{PO2} M.~A.~Olshanetsky and A.~M.~Perelomov.  Quantum systems
  related to root systems, and radial parts of Laplace operators.
  \emph{Funct. Anal. Appl. } \textbf{12} (1978), 121--128.

\bibitem{pakhomov} V.~F. Pahomov.  \newblock Automorphisms of the
  tensor product of an abelian and a {G}rassmann algebra.  \newblock
  {\em Math. Notes} \textbf{16} (1974), 624--629. Transl. from:
  \emph{Matem. Zametki} \textbf{16} (1) (1974), 65--74.

\bibitem{ramond:dual71} P.~Ramond.  \newblock Dual theory for free
  fermions.  \newblock {\em Phys. Rev. D (3)} \textbf{3} (1971),
  2415--2418.

\bibitem{salamstrathdee:supergauge} A.~Salam and J.~Strathdee.
  \newblock Super-gauge transformations.  \newblock {\em Nuclear
    Phys. B} \textbf{76} (1974), 477--482.


\bibitem{Sch} M.~Schlichenmaier. Berezin-Toeplitz quantization of
  compact K\"ahler manifolds. In: Quantization, Coherent States and
  Poisson Structures, Proc. XIVth Workshop on Geometric Methods in
  Physics (Bia{\l}owie\.{z}a, Poland, 9--15 July 1995),
  A. Strasburger, S. T. Ali, J.-P. Antoine, J.-P. Gazeau, and
  A. Odzijewicz, eds., Polish Scientific Publisher PWN (1998),
  101--115.

\bibitem{schmitt:ident} Th.~Schmitt.  \newblock Some identities for
  {B}erezin's function.  \newblock In {\em Seminar Analysis, 1981/82},
  pages 146--161. Akad. Wiss. DDR, Berlin, 1982.


\bibitem{schwinger:dynamprinc53} J.~Schwinger.  \newblock A note on
  the quantum dynamical principle.  \newblock {\em Philos. Mag. (7)},
  \textbf{44} (1953), 1171--1179.

\bibitem{schwinger:grassmann62} J.~Schwinger.  \newblock Exterior
  algebra and the action principle, I.  \newblock {\em PNAS},
  \textbf{48}(4) (1962), 603--611.

\bibitem{Seg2} I.~E.~Segal. Tensor algebras over Hilbert
  spaces. I. \emph{Trans. Amer. Math. Soc.} \textbf{81} (1956),
  106--134.


\bibitem{Seg} I.~E.~Segal. Foundations of the theory of dynamical
  systems of infinitely many degrees of freedom. I.  \emph{
    Mat.-Fys. Medd. Danske Vid. Selsk.}  \textbf{31}, no. 12, (1959)
  39 pp.


\bibitem{Sek} J.~Sekiguchi. Zonal spherical functions on some
  symmetric spaces. \emph{Publ. Res. Inst. Math. Sci., Kyoto Univ.}
  \textbf{12}, Suppl., (1977), 455--459.

\bibitem{Ser} A.~N.~Sergeev. Calogero operator and Lie superalgebras.
  \emph{Theoretical and Mathematical Physics} \textbf{131} (3) (2002),
  747--764

\bibitem{Sha} D.~Shale. Linear symmetries of free boson
  fields. \emph{Trans. Amer. Math. Soc.}  \textbf{103} (1962),
  149--167.

\bibitem{ShSt} D.~Shale and W.~F.~Stinespring. Spinor representations
  of infinite orthogonal groups.  \emph{J. Math. Mech. } \textbf{14}
  (1965), 315--322.

\bibitem{SSV} D. Shklyarov, S. Sinel'shchikov, and L. Vaksman. \newblock A q-analogue of the Berezin quantization method, \newblock \emph{Lett. Math. Phys.} \textbf{49}(3) (1999), 253--261.

\bibitem{Soe} W.~Soergel. An irreducible not admissible Banach
  representation of $SL(2,\mathbb R)$. \emph{Proc. Amer. Math. Soc.}
  \textbf{104}, No.4 (1988), 1322--1324.

\bibitem{UU} A.~Unterberger, H.~Upmeier.  The Berezin transform and
  invariant differential operators.  \emph{Comm. Math. Phys.}
  \textbf{164} (1994), 563--597.



\bibitem{volkov:akulov72} D.V. Volkov and V.P. Akulov.  \newblock
  Possible universal neutrino interaction.  \newblock {\em JETP
    Letters} \textbf{16}(11) (1972), 438--440.


\bibitem{volkov:akulov74} D.V. Volkov and V.P. Akulov.  \newblock
  Goldstone fields with spin 1/2.  \newblock {\em Theoretical and
    Mathematical Physics} \textbf{18}(1) (1974):28--35.

\bibitem{volkov:soroka74} D.~V. Volkov and V.~A. Soroka.  \newblock
  Gauge fields for a symmetry group with spinor parameters.  \newblock
  {\em Theoretical and Mathematical Physics} \textbf{20}(3) (1974),
  829--834.

\bibitem{tv:git} Th. Voronov.  \newblock {\em Geometric Integration
    Theory on Supermanifolds}, volume~9 of {\em
    Sov. Sci. Rev. C. Math. Phys.}  \newblock Harwood Academic Publ.,
  1992.

\bibitem{wesszumino:supergauge74} J.~Wess and B.~Zumino.  \newblock
  Supergauge transformations in four dimensions.  \newblock {\em
    Nuclear Phys. B} \textbf{70}:39--50, 1974.

\bibitem{wesszumino:lagrmodel74} J.~Wess and B.~Zumino.  \newblock A
  {Lagrangian} model invariant under supergauge transformations.
  \newblock {\em Phys. Lett.} \textbf{49B}(1) (1974), 52--54.



\bibitem{Wei} A.~Weil. Sur certains groupes d'operateurs unitaires.
  \emph{Acta Math. } \textbf{111} (1964), 143--211.


\bibitem{Zhe} D.~P.~Zhelobenko.  Operational calculus on a semisimple
  complex Lie group.  \emph{Izv. Akad. Nauk SSSR Ser. Mat.}
  \textbf{33} (1969), 931--973; English translation in
  \emph{Mat.-USSR-Izvestia} \textbf{3} (1971), 881--916.

\bibitem{ZhN} D.~P.~Zhelobenko and M.~A.~Naimark. A characterization
  of completely irreducible representations of a semisimple complex
  Lie group. \emph{Sov. Math.  Dokl.} \emph{7} (1966), 1403--1406 ;
  translation from \emph{Dokl. Akad. Nauk SSSR} \textbf{171} (1966),
  25--28.
\end{thebibliography}

% ------------------------------------------------------------------------
\end{document}